\documentclass[10pt,oneside]{amsart}

\usepackage{amscd}
\usepackage{amssymb}
\usepackage[latin1]{inputenc}   
\usepackage[pdftex]{graphicx}   

\def\Z{\mathbb{Z}}
\def\Q{\mathbb{Q}}
\def\ZZ{\Z \oplus \Z}
\def\R{\mathbb{R}}
\def\C{\mathbb{C}}
\def\H{\mathbb{H}}
\def\E{\mathbb{E}}
\def\SS{\mathbb{S}}
\def\KB{\mathbb{KB}}

\def\a{\alpha}          
\def\b{\beta}
\def\g{\gamma}
\def\d{\delta}
\def\e{\varepsilon}
\def\f{\phi}
\def\vf{\varphi}
\def\l{\lambda}

\def\r{\rho}
\def\s{\sigma}
\def\w{\omega}
\def\G{\Gamma}

\def\P{\pi_{1}}         
\def\cal{\mathcal}
\def\T{\cal{T}}
\def\ol{\overline}
\def\ul{\underline}
\def\lb{\lbrack}
\def\rb{\rbrack}
\def\wt{\widetilde}
\def\bs{$\square$}
\def\o{\mathbf{o}}
\def\ex{\mathbf{e}}
\def\lgr{\mathrm{lg}}

\newtheorem{thm}{Theorem}[section]
\newtheorem{cor}{Corollary}[section]
\newtheorem{lem}{Lemma}[section]
\newtheorem{prop}{Proposition}[section]

\title{Conjugacy problem in groups of oriented geometrizable
3-manifolds}

\begin{document}

\maketitle

\begin{center}
{\sc Jean-Philippe PR\' EAUX}\footnote[1]{{ Jean-Philippe
PR\'EAUX}, \\
-- Ecole de l'air,\ CREA,\ 13661 Salon de Provence air, France\\
-- Centre de Math\'ematiques et d'informatique, Universit\'e de
Provence, 39 rue F.Joliot-Curie, 13453 marseille cedex 13, France.
\ {\it E-mail :} \ preaux@cmi.univ-mrs.fr}
\end{center}

\begin{abstract}
The aim of this paper is to show that the fundamental group of an
oriented 3-manifold which satisfies Thurston's geometrization
conjecture, has a solvable conjugacy problem. In other terms, for
any  such 3-manifold $M$, there exists an algorithm which
can decide for any couple of elements $u,v$ of $\P(M)$  whether
 $u$ and $v$ are in the same conjugacy class of $\P(M)$ or not. More
topologically, the algorithm decides for any two loops in $M$, whether they are freely
homotopic or not.
\end{abstract}


\section*{\bf Introduction}

Since the work of M.Dehn (\cite{dehn1}, \cite{dehn2},
\cite{dehn3}), the Dehn problems (more specifically the word
problem and the conjugacy problem) have taken a major importance
in the developments of combinatorial group theory
 all along the $\mathrm{XX}^{th}$ century. But while elements of this theory provide
 straightforward solutions in special cases, there still remains
 important classes of groups for which actual techniques fail to
 give a definitive answer. This is particularly true, when one
 considers the conjugacy problem. This last problem seems to be
 much more complicated  than the word problem, and even
  very simple cases are still open or admit a negative answer.

  Let us focus on fundamental groups of compact manifolds, case
  which motivated the introduction by Dehn of these problems. M.
  Dehn has solved the word and conjugacy problems for groups of
  surfaces. It is well known that these two problems are in
  general insolvable in groups of $n$-manifolds for $n\geq 4$
  (as a consequence of the facts that each finitely
  presented group is the fundamental group of some $n$-manifold
  for $n\geq 4$ fixed, and of the general insolvability of these
  two
  problems for an arbitrary f. p. group). In case of dimension 3,
  not all f. p. groups occur as fundamental groups of 3-manifolds, but the problems are still
  open for this class, despite many improvements. To solve the word problem,
  one needs to suppose a further condition, namely Thurston's
  geometrization conjecture. In such a case, the automatic group
  theory gives a solution to the word problem, but fails to solve
  the conjugacy problem (cf. \cite{epstein}). Some special cases admit nevertheless a
  solution. Small cancellation theory provides a solution for alternating links, and
    biautomatic group theory provides solutions for hyperbolic manifolds and almost all Seifert fiber spaces ; the   best result we know applies to irreducible 3-manifolds with non-empty boundary : they admit a locally $CAT(0)$ metric and hence their groups have solvable conjugacy problem (\cite{bridson}).
  A major improvement has been the solution of Z.Sela (\cite{sela}) for
  knot groups. In his paper, Sela conjectures "the method seems to
  apply to all 3-manifolds satisfying Thurston's geometrization
  conjecture". We have been inspired by his work, to show the main
  result of this paper :\smallskip\\
  {\bf Main Theorem.}   {\sl
  The group of an oriented 3-manifold satisfying Thurston's
  geometrization conjecture has a solvable conjugacy problem.
  }\smallskip\\
 \indent Let us first give precise definitions of the concept involved.
  Let $G=<X|R>$ be a group given by a finite presentation. The
  {\it word problem}, consists in finding an algorithm which decides for
  any couple of words $u,v$ on the generators, if $u=v$ in $G$.
  The {\it conjugacy problem}, consists in finding an algorithm which decides for
  any couple of words $u,v$ on the generators, if $u$ and $v$ are conjugate in $G$ (we shall write $u\sim v$),
  that is if there exists $h\in G$,
  such that $u=h.v.h^{-1}$ in $G$. Such an algorithm is
  called a {\it solution} to the corresponding problem. It turns
  out that the existence of a solution to any of these problems,
  does not depend upon the finite presentation of $G$ involved.
  Novikov (\cite{novikov},  1956)
  has shown that a solution to the word problem does not in general
  exist. Since a solution to the conjugacy problem provides a
  solution to the word problem (to decide if $u=v$ just decide if
  $uv^{-1}\sim 1$), the same conclusion applies to the conjugacy
  problem. Moreover, one can construct many examples of groups
  admitting a solution to the word problem, and no solution to the
  conjugacy problem (cf. \cite{miller}).

  When one restricts to the fundamental group of a manifold,
  solving the word problem is equivalent to deciding for any
  couple of based point loops, whether they are or not
  homotopic with base point fixed, while solving the conjugacy
  problem is equivalent to deciding whether two loops are freely
  homotopic.

By a 3-{\it manifold} we mean a compact connected
$\ul{\mathrm{oriented}}$ manifold with boundary, of dimension 3.
According to the Moise theorem (\cite{moise}), we may use
indifferently PL, smooth, or topological locally smooth, structures on such a
3-manifold.

A 3-manifold is said to satisfy Thurston's geometrization property
(we will say that the manifold is {\it geometrizable}), if it
decomposes in its canonical topological decomposition ---along
essential spheres (Kneser-Milnor), disks (Dehn lemma) and tori (Jaco-Shalen-Johanson)---
 into pieces whose interiors admit
a riemanian metric complete and locally homogeneous. (In the following we shall speak improperly of a geometry on a
3-manifold rather than on its interior.)

Thurston's has conjectured that all 3-manifolds are geometrizable.
This hypothesis is necessary to our work, because otherwise one
can at this date classify neither 3-manifolds nor their
groups.

  To solve the conjugacy problem, we will first use the classical topological decomposition, as well as the
classification of geometrizable 3-manifolds, to reduce the conjugacy
problem to the restricted case of closed irreducible 3-manifolds,
which are either Haken (i.e. irreducible and containing a properly
embedded 2-sided incompressible surface, cf.
\cite{jaco},\cite{js}), a Seifert fibered space (\cite{seifert},
\cite{jaco}), or modelled on $SOL$ geometry (cf. \cite{scott}).
That is, we show that if the group of any 3-manifold lying in such
classes has a solvable conjugacy problem, then the same conclusion
applies to any geometrizable 3-manifold (lemma \ref{red}). The cases
of Seifert fibered spaces and $SOL$ geometry are rather easy, and
we will only sketch solutions respectively in \S 5.3 and \S 7 ;
the inquiring reader can find detailed solutions in my PhD thesis,
\cite{thesis}. The Haken case constitutes the essential
difficulty, and will be treated in details. We can further suppose
that the manifold is not a torus bundle, because in such a case
the manifold admits either a Seifert fibration or a geometric
structure modelled on $SOL$ geometry. As explained above, we solve
the conjugacy problem in the group of a Haken closed manifold by
essentially applying the strategy used by Z.Sela to solve the case
of knot groups.

Let's now focus on the  case of a Haken closed manifold $M$. The
JSJ decomposition theorem asserts that there exists a family of
essential tori $W$, unique up to ambient isotopy, such that if one
cuts $M$ along $W$, one obtains pieces which are either Seifert
fibered spaces or do not contain an essential torus (we say the
manifold is atoroidal). According to Thurston's geometrization
theorem (\cite{thu}), each atoroidal piece admits a hyperbolic
structure with finite volume. This decomposition of $M$ provides a
decomposition of $\P(M)$ as a fundamental group of a graph of
groups, whose vertex groups are the $\P$ of the pieces obtained,
and edge groups are free abelian of rank 2. We then establish a
conjugacy theorem (in the spirit of the conjugacy theorem for
amalgams or HNN extensions), which characterizes conjugate
elements in $\P(M)$ (theorem \ref{conjthm}). This result, together
with the algebraic interpretation of the lack of essential annuli
in the pieces obtained (proposition \ref{annuli}), and with a
solution to the word problem, allows us to reduce the conjugacy
problem in $\P(M)$ to three algorithmic problems in the groups of
the pieces obtained : the conjugacy problem, the boundary
parallelism problem, and the 2-cosets problem (theorem
\ref{reduce}). Suppose $M_1$ is a piece, and $T$ is a boundary
subgroup of $\P(M_1)$ (that is $T=i_*(\ZZ)\subset\P(M_1)$ for some
embedding $i:S^1\times S^1\longrightarrow \partial M_1$). The
{\sl boundary parallelism problem in} $\P(M_1)$ consists in finding for
any element $\w\in \P(M_1)$ all the elements of $T$ conjugate to
$\w$ in $\P(M_1)$. Suppose $T_1,T_2$ are two boundary subgroups
(possibly identical), the {\sl 2-cosets problem} consists in finding for
any couples $\w,\w'\in \P(M_1)$, all the elements $c_1\in T_1$,
$c_2\in T_2$, such that $\w=c_1.\w'.c_2$ in $\P(M_1)$. We then
solve these algorithmic problems, separately, in the Seifert case,
and in the hyperbolic case, providing a solution to the conjugacy
problem in $\P(M)$.

In the Seifert case a solution can be easily established, by using
the existence of a normal infinite cyclic subgroup $N\subset
\P(M_1)$, such that $\P(M_1)/N$ is virtually a surface group.
Algorithms in $\P(M_1)$ can be reduced to similar algorithms in
$\P(M_1)/N$, providing solutions (propositions \ref{bound_seif}
and \ref{2coset_seif}). A solution to the conjugacy problem in the
Seifert case, already solved in almost all cases by biautomatic
group theory and presenting no difficulty in the (few) remaining
cases (namely $NIL$ geometry), will only be sketched in \S 5.3.

In the case of a hyperbolic piece $M_1$,  biautomatic group
theory solves the conjugacy problem. The 2-cosets problem will be
solved using the hyperbolic structure of $\P(M_1)$ relative to its
boundary subgroups (in the sense of Farb, \cite{farb}). To solve
the boundary parallelism problem, we will make use of Thurston's
hyperbolic surgery theorem to obtain two closed hyperbolic
manifolds by Dehn filling on $M_1$, and then reducing the problem
in $\P(M_1)$ to analogous problems in the groups of these two
manifolds. Solutions will then be provided using word hyperbolic
group theory.

\section{\bf Reducing the problem}

The aim of this section is to reduce the conjugacy problem in the
group of a geometrizable (oriented) 3-manifold to the same problem in
the more restricted case of a closed irreducible 3-manifold which
is either Haken, or a Seifert fibered space, or modelled on $SOL$
geometry. That is, if the conjugacy problem is solvable in the
group of any such 3-manifold, then it is also solvable in the group of any geometrizable 3-manifold. This result
constitutes lemma \ref{red} ; the reduction will be done in three
steps : first reducing to a closed manifold by "doubling the
manifold", then to an irreducible closed manifold by using the
Kneser-Milnor decomposition, to finally conclude with the
classification theorem for closed irreducible geometrizable
3-manifolds.

\subsection{Reducing to the case of a closed manifold}
Suppose $M$ is a 3-manifold with non-empty boundary. Consider an
homeomorphic copy $M'$ of $M$, and an homeomorphism $\vf : M
\longrightarrow M'$. Glue $M$ and $M'$ together along the
homeomorphisms induced by $\vf$ on the boundary components, to
obtain a closed 3-manifold, which  will be called $2M$.  We can
reduce the conjugacy problem in $\P(M)$ to the conjugacy problem
in $\P(2M)$, by using the following lemma.

\begin{lem} $\P(M)$ naturally embeds in
$\P(2M)$. Moreover, two elements $u,v\in \P(M)$ are conjugate in
$\P(M)$ if and only if they are conjugate in $\P(2M)$.
\end{lem}

\noindent {\bf Proof.} For more convenience $2M$ can be seen as
$2M=M\cup M'$, with $\partial M=\partial M'=M\cap M'$. There
exists a natural homeomorphism $\f : M\longrightarrow M'$ which
restricts to identity on $\partial M$. Consider the natural (continuous)
map $\pi :2M\longrightarrow M$ defined by $\pi_{|M}=Id_M$ and $\pi_{|M'}=\f^{-1}$.

The inclusion $M\subset 2M$ induces a group homomorphism
$\P(M)\longrightarrow\P(2M)$ . We need first to prove that this
map is injective. It suffices to show that any loop in $M$,
contractile in $2M$, is in fact contractile in $M$. Consider a
loop $l$
 in $M$, that is $l:S^1\longrightarrow M$, and suppose that there
exists a map $h:D^2\longrightarrow 2M$ such that $h$ restricted to
$S^1=\partial D^2$ is $l$ ; note $\ol{h}=\pi\circ h$.
Since the loop $l$ lies
in $M$, $\ol{h}$ restricted to $\partial D^2$ is $l$. Hence $l$ is
contractile in $M$, which proves the first assertion.

We now prove the second assertion. Since $\P(M)$ is a subgroup of
$\P(2M)$, the direct implication is obvious. we need to prove the
converse. Consider two loops $l_u$ and $l_v$ in $M$ which
represent respectively $u$ and $v$.  Suppose also, that $u$ and
$v$ are conjugate in $\P(2M)$. So $l_u$ and $l_v$ are freely
homotopic in $2M$. Thus, there exists a map $f: S^1\times I
\longrightarrow 2M$, such that $f$ restricted to $S^1\times 0$ is
$l_u$ and $f$ restricted to $S^1\times 1$ is $l_v$ ; note
$\ol{f}=\pi\circ f$. Since $l_u,l_v$ lie in $M$, $\ol{f}$ is an homotopy from
$l_u$ to $l_v$ in $M$, and so $u$ and $v$ are conjugate in $\P(M)$.\hfill\bs\\

Suppose one needs to solve the conjugacy problem in $\P(M)$ where
$M$ is 3-manifold with non-empty boundary. By doubling the
manifold $M$ along its boundary, one obtains the closed 3-manifold
$2M$. If the conjugacy problem $\P(2M)$ admits a solution, then
one can deduce a solution in $\P(M)$. Consider $u$ and $v$ in
$\P(M)$. Under the natural embedding $\P(M) \hookrightarrow
\P(2M)$, $u$ and $v$ can be seen as elements of $\P(2M)$. With the
preceding lemma, one only needs to check if $u$ and $v$ are
conjugate in $\P(2M)$ to determine if they are conjugate in
$\P(M)$. Hence conjugacy problem in $\P(M)$ reduces to conjugacy
problem in $\P(2M)$. Together with the following lemma, the
conjugacy problem in geometrizable (oriented) 3-manifolds reduces
to the conjugacy problem in closed (oriented) geometrizable
3-manifolds.

\begin{lem}
 If $M$ is geometrizable, then so is $2M$.
\end{lem}

\noindent {\bf Proof.} A 3-manifold is geometrizable if and only if all
of its prime factors are geometrizable. Together with Kneser-Milnor theorem,
if $M$ splits as $M=\# M_i$, then $M$ is geometrizable if and only if each
of its (non necessarily prime) factors $M_i$ are geometrizable ;
 this observation will be denoted by $(*)$.

We suppose $M$ to be geometrizable, and want to show that $2M$ is
geometrizable. We note $C_i$ (resp.  $D_i$) the prime factors of
$M$ with empty (resp. $\not=\emptyset$) boundary ; then $M=(\#
C_i)\# (\# D_j)$ and $2M= (\# C_i)\# (\# C'_i) \# (\# 2D_j)$ where
$C'_i$ are homeomorphic copies of $C_i$ ; be careful that the
$2D_j$ are not necessarily prime. Using $(*)$ it suffices to show
that all the $2D_j$ are geometrizable. Hence we will suppose in
the following that $M$ is prime with non-empty boundary ; we
suppose besides that $\partial M \not\supset S^2$ cause otherwise
$M=B^3$ and $2M=S^3$.

If $M$ is not $\partial$-irreducible then it must contain an essential disk.
 If $M$ contains a separating essential disk, $M=M_1\#_{D^2}M_2$, then $2M$
  contains a separating essential sphere and splits non trivially as $2M=2M_1\# 2M_2$.
  Hence, using $(*)$ and the fact that a 3-manifold is not infinitely decomposable
  as a non trivial connected sum,
we will suppose that $M$ does not contain any separating essential disk. In particulary $2M$
is prime, and if $M$ would contain a non separating disk then
$2M$ would be a sphere bundle over $S^1$ ;
so that we  moreover suppose $M$ to be $\partial$-irreducible. Under these hypothesis
$2M$ is Haken, and according to Thurston's geometrization theorem, is geometrizable.
\hfill $\square$

\subsection{Reducing to the case of an irreducible manifold}

Suppose now that $M$ is a closed 3-manifold. According to the
Kneser-Milnor theorem (cf. \cite{hempel}), $M$ admits a unique
decomposition as a connected sum of prime manifolds,
$M=M_1\#M_2\#\cdots \#M_n$, where each $M_i$ is either
irreducible, or homeomorphic to $S^2\times S^1$. Its fundamental
group $\P(M)$ decomposes in a free product of the $\P(M_i)$ for
$i=1,2,\ldots ,n$, that is $\P(M)=\P(M_1)\ast\P(M_2)\ast\cdots
\ast\P(M_n)$.

Applying the conjugacy theorem for a free product (cf.
\cite{mks}), $\P(M)$ has a solvable conjugacy problem, if and only
if each of the $\P(M_i)$ has a solvable conjugacy problem.
Moreover $\P(S^2\times S^1)$ is infinite cyclic and therefore
admits a solution to the conjugacy problem. So the conjugacy
problem in a (closed, geometrizable) 3-manifold group reduces to the
conjugacy problem in all (closed, geometrizable) irreducible
3-manifolds groups. Together with the last lemma, one obtains

\begin{lem} The conjugacy problem in groups of oriented (resp. geometrizable) 3-manifolds,
reduces to conjugacy problem in groups of closed and irreducible
oriented (resp. geometrizable) 3-manifolds.
\end{lem}

\subsection{Reducing to particular 3-manifolds}

The final step in this reduction, comes from the classification
theorem for closed irreducible geometrizable manifolds together with
all the already known results on the conjugacy problem in the
groups of such manifolds.
The following result constitutes theorem 5.3 of \cite{scott}
slightly adapted to the oriented case.

\begin{thm}
\label{scotty}
 Let $M$ be an irreducible, closed geometrizable
(oriented) 3-manifold. Then $M$ satisfies one
of the following conditions : \medskip\\
\indent (i) $M$ is Haken.\medskip\\
\indent (ii) $M$ is hyperbolic.\medskip\\
\indent (iii) $M$ is modelled on $SOL$ geometry. This happens
exactly when $M$ is finitely covered by a $S^1\times S^1$-bundle
over $S^1$, with hyperbolic gluing map.  In particular, either $M$
is itself a $S^1\times S^1$-bundle over $S^1$, or $M$ is the union
of two twisted $I$-bundles over the Klein bottle. In this case,
$M$ is Haken.\medskip\\
\indent (iv) $M$ is modelled on $\SS^3$, $\E^3$, $\SS^2\times
\E^1$, $\H^2\times \E^1$, $NIL$, or on the universal cover of
$SL(2,\R)$. This happens exactly when $M$ is a Seifert fibered
space. In this case $M$ is a $S^1$-bundle with base an orbifold
$O_2$, and if  $e$ refers to the Euler number of the bundle, and
$\chi$ to the Euler characteristic of the base $O_2$, then the
geometry of $M$ is characterized by $e$ and $\chi$, following the
table below :
\begin{center}
\begin{tabular}{cccc}
&  $\qquad \chi>0\qquad $ & $\qquad\chi=0\qquad$ &
$\qquad\chi<0\qquad$
\\
&&&\\
$e=0\quad$ & $\SS^2\times \E^1$ &
$\E^3$ & $\H^2\times\E^1$\medskip\\
$e\not= 0\quad$  & $\SS^3$ & $NIL$ & $\wt{SL}(2,\R)$\\
\end{tabular}
\end{center}
\end{thm}

Note that conditions $(ii)$, $(iii)$ and $(iv)$ are disjoint
(according to the unicity of the geometry involved for a
particular 3-manifold) while condition $(i)$ is not disjoint from
conditions $(ii)$, $(iii)$ and $(iv)$. For example, obviously a
torus bundle $M$ over $S^1$ is Haken, while $M$ is modelled on
$SOL$ when its gluing map is Anosov, or $M$ is a Seifert
fibered space (and so satisfies condition $(iv)$) in the case of a
reducible or periodic
gluing map.\\

In the case of a closed hyperbolic manifold $M$, $M$ is
the orbit space of a cocompact action of a discrete subgroup of
$PSL(2,\C)$ on $\H^3$. So, $\P(M)$ is word-hyperbolic (in the
sense of Gromov, cf. \cite{gro},\cite{cdp}), and admits a (very
efficient) solution to the conjugacy problem. So we only need to
solve the problem in the remaining cases $(i)$, $(iii)$ and
$(iv)$. We have finally obtained the main result of this section :
\begin{lem}
\label{red} The conjugacy problem in groups of oriented geometrizable
3-manifolds, reduces to groups of oriented closed 3-manifolds
which are either Haken, or a Seifert fibered space, or modelled on
$SOL$ geometry.
\end{lem}

In the case of a Seifert fiber space $M$, its group $\P(M)$ is
biautomatic, unless $M$ is modelled on $NIL$ (cf. \cite{nr1},
\cite{nr2}). Hence  biautomatic group theory provides a solution
to the conjugacy problem in almost all cases. The remaining cases
(those modelled on $NIL$) can easily be solved by direct methods.
We will only sketch a solution in \S 5.3.

In the case of a 3-manifold $M$ modelled on $SOL$ geometry, $M$ is
either a $S^1\times S^1$-bundle over $S^1$ or obtained by gluing
two twisted $I$-bundles over $\KB_2$ along their boundary ; in
particular $M$ is Haken. We distinguish the $SOL$ case from the
Haken case, cause we solve separately what we shall call the {\it
generic Haken case} (a Haken closed manifold which is neither a
$S^1\times S^1$-bundle over $S^1$, nor obtained by gluing two
$I$-twisted bundles over $\KB_2$ along their boundary), from the
remaining Haken cases (namely a $S^1\times S^1$-bundle over $S^1$,
or two $I$-twisted bundles over $\KB_2$ glued along their
boundary).

The reason for such a distinction, is that in this last non generic case our general strategy
fails (because the JSJ decompositions may not be $k$-acylindrical for some $k>0$, see \S 4.3). Nevertheless
a
solution can easily be established  : either they are Seifert fibered, or the
conjugacy problem reduces easily to solving elementary equations in
$SL(2,\Z)$. We will only sketch a solution in \S 7.

The main part of our work will be devoted to the Haken generic
case, which represents the major difficulty, and will be treated
in details. While solutions in the remaining cases are only
sketched in \S 5.3 and $\S 7$, the inquiring reader can find
detailed solutions in my PhD. thesis \cite{thesis}, \S 5.5 and \S 7.1.

\section{\bf The group of a Haken closed manifold}

We study in this section the fundamental group of a Haken closed
manifold $M$. We see how the JSJ decomposition of $M$ provides a
splitting of $\P(M)$ as a fundamental group of a graph of groups.

\subsection{JSJ Decomposition.}
Let $M$ be a Haken closed manifold. The JSJ  theorem (cf.
\cite{js}) asserts that there  exists an essential surface $W$
embedded in $M$ (possibly $W=\emptyset$), whose connected
components consist of tori, such that if one cuts $M$ along $W$,
one obtains a 3-manifold, whose connected components -- called the
{\it pieces} --  are either Seifert fibered spaces, or atoroidal
(i.e. do not contain an essential torus). Moreover, $W$ is minimal
up to ambient isotopy, in the class of surfaces satisfying the above
conditions. The manifold $M$ can then be reconstructed by gluing
the pieces along their boundary components.

The minimality of $W$ has two consequences which will be essential
in the following of our work. First if $M_1$, $M_2$ are two
Seifert pieces glued along one boundary component (in order to
reconstruct $M$), then the gluing map sends a regular fiber of
$M_1$ onto a loop of $M_2$ which cannot be homotoped to any
regular fiber in $M_2$ : cause otherwise one can extend the
fibration on the gluing of $M_1$ and $M_2$, contradicting the
minimality of the decomposition. This fact is resumed in the
following lemma :

\begin{lem}
\label{fibers}
 Suppose $M_1$ and $M_2$ are two Seifert pieces in the JSJ decomposition of $M$, which are
 glued along one boundary component. Then the gluing map sends a regular fiber of $M_1$ to a loop
 which cannot be homotoped in $M_2$ to a regular fiber.
\end{lem}

The second consequence excludes in almost all cases, pieces
homeomorphic to a thickened torus. The proof is immediate since
gluing one piece $N$ with $S^1\times S^1\times I$ along one
boundary component does not change the homeomorphism class of $N$
while it increases the number of components of $W$.

\begin{lem}
\label{s1s1i} Suppose $M$ is a Haken closed manifold which is not
homeomorphic to a $S^1\times S^1$-bundle over $S^1$. Then none of
the pieces of the JSJ decomposition of $M$ is homeomorphic to
$S^1\times S^1\times I$.
\end{lem}

We now introduce some notations. Suppose $W$ is non-empty, and
admits as connected components $\cal{T}_1$, $\cal{T}_2,\ldots
,\cal{T}_q$. since $W$ is two-sided in $M$, each of the
$\cal{T}_i$ admits a regular neighborhood in $M$ homeomorphic to a
thickened torus, which we shall note $N(\cal{T}_i)$, chosen in
such a way that all the $N(\cal{T}_i)$ (for $i=1,2,\ldots ,q$) are
disjoints one to each other. Now when we say "cutting $M$ along
$W$" we mean considering the compact 3-manifold $\s_W(M)$ defined
as :
$$\s_W(M)=M - \bigcup int(N(\cal{T}_i))$$
The homeomorphism class of $\s_W(M)$ does not depend on the
neighborhoods involved. The connected components of $M$ have
non-empty boundary when $W\not=\emptyset$. We shall call them the
{\it pieces} of the decomposition, and  name them as : $M_1, M_2,
\ldots ,M_n$.

According to Thurston's geometrization theorem (cf. \cite{thu}),
the atoro\" \i dal pieces admit a hyperbolic structure with
finite volume. We shall call them the {\it hyperbolic pieces}.

There exists a canonical map associated to the JSJ decomposition,
called the {\it identification map}, $r:\s_W(M)\longrightarrow M$,
which is such that $r$ restricted to $int(\s_W(m))$ is an
homeomorphism, and each preimage $r^{-1}(\cal{T}_i)$ of
$\cal{T}_i$ consists of two homeomorphic copies of $S^1\times S^1$
in $\partial \s_W(M)$, which we shall (arbitrarily) call
$\cal{T}_i^-$ and $\cal{T}_i^+$.

Then $M$ can be  reconstructed from $\s_W(M)$. There exists two
homeomorphisms :
\begin{gather*}
\r_i^-:S^1\times S^1\longrightarrow \cal{T}_i^-\\
\r_i^+:S^1\times S^1\longrightarrow \cal{T}_i^+
\end{gather*}
such that the following diagram commutes :
$$
\begin{CD}
S^1\times S^1 @>\r_i^->> \cal{T}_i^-\\
@V\r_i^+VV  @VVrV \\
\cal{T}_i^+ @>>r> \cal{T}_i
\end{CD}
$$
Define $\r_i:\cal{T}_i^-\longrightarrow \cal{T}_i^+$ by
$\r_i=\r_i^+\circ (\r_i^-)^{-1}$, which will be called {\it the
gluing map associated to the torus} $\cal{T}_i$. Then $M$ is
homeomorphic to the manifold obtained by gluing $\s_W(M)$ on its
boundary, according to the gluing maps $\r_1,\r_2,\ldots ,\r_q$.

\subsection{Splitting the fundamental group}

The JSJ decomposition of $M$ provides a splitting of $\pi_1(M)$ as a fundamental group of a graph of groups.
The information needed to define the graph of groups comes
directly from the information characterizing the JSJ decomposition
of $M$, namely the pieces obtained in this decomposition, and the
associated gluing maps.

A JSJ decomposition
$W=\T_1\cup\T_2\cup \cdots \cup \T_q$  of $M$, and  the pieces $M_1,M_2,\ldots
,M_n$ obtained in this decomposition naturally provide a splitting of $M$ as a graph of space
$(\cal{M},X)$ (cf. \cite{scott-wall}).
The underlying graph $X$ has one vertex $v(M_i)$ (which shall also be noted $v_i$) for
each piece $M_i$ ; given a vertex $v$ of $X$ we shall denote by $M(v)$ the associated piece.
To each connected component $\T_j$ of $W$ correspond two  edges of $X$
$e(\T_j)$ and $e(\T_j)^{-1}$ inverse one to the other (we shall also note respectively $e_j$
and $\bar{e}_j$), and if
$M_k,M_l$ are such that $\T_j^-\subset\partial M_k$
and $\T_j^+\subset\partial M_l$, then $e(\T_j)$ has origin
$\o(e(\T_j))=v(M_k)$ and extremity $\ex(e(\T_j))=v(M_l)$. The gluing maps associated to the
edge $e(\T_j)$ are given respectively by $\r_i^-$ and $\r_i^+$.
\begin{figure}[ht]
\center{\includegraphics[scale=0.7]{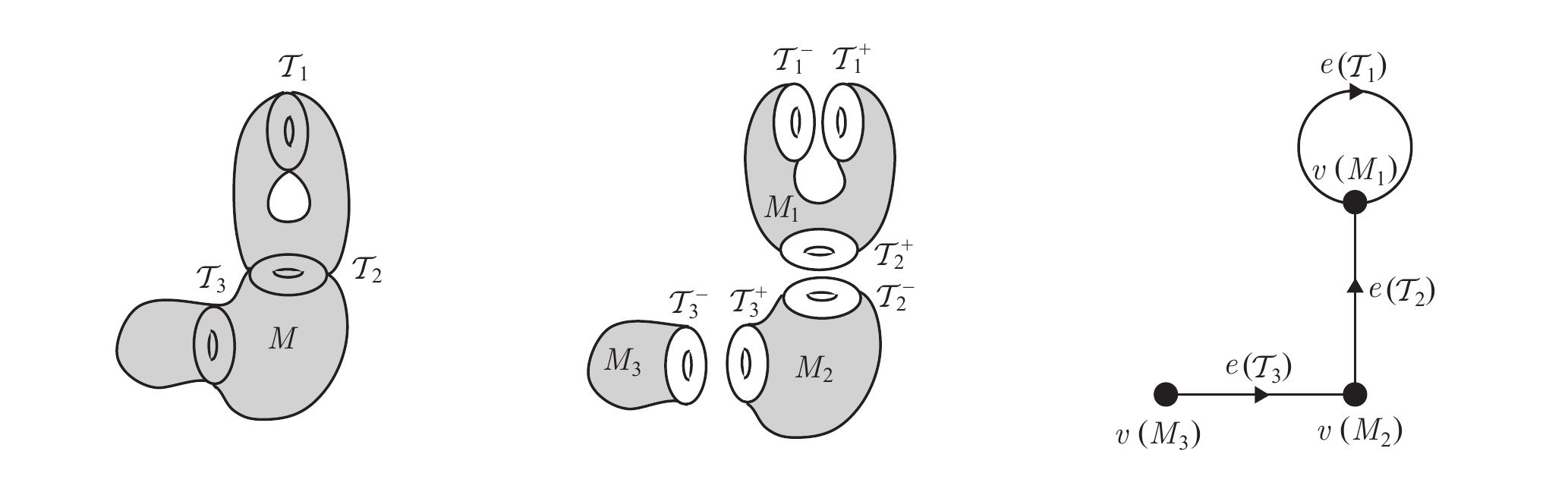}}
\caption{Construction of the graph $X$}
\end{figure}
One naturally obtains the graph of groups $(\cal{G},X)$ by assigning to each vertex $v_i=v(M_i)$
the group $\P(M_i,*_i)$ for some $*_i\in M_i$ --- which
we shall denote by $G(v_i)$  and call a {\it
vertex group} --- and to each edge $e_j$ the group $G(e_j)$ isomorphic to
$\Z\oplus \Z$ --- called an {\it edge group}. For each edge $e_i$, say with $\o(e_i)=v_0$ and
$\ex(e_i)=v_1$, fix one base point $\ul{*}_i^- \in \T_i^-$ (resp. $\ul{*}_i^+ \in \T_i^+$ such
that $r(\ul{*}_i^-)=r(\ul{*}_i^+)=\ul{*}_i\in\T_i$) and a path from $*_0$ to $\ul{*}_i^-$
in $M(v_0)$ (resp. from $*_1$ to $\ul{*}_i^+$ in $M(v_1)$) ; it defines
two monomorphisms $\vf_j^-:G(e_j)\hookrightarrow G(v_0)$
and $\vf_j^+:G(e_j)\hookrightarrow G(v_1)$ by
$\vf_j^-=(\r_j^-)_\ast$, $\vf_j^+=(\r_j^+)_\ast$, which shall also be denoted by $\vf_{e_j}^-$,
$\vf_{e_j}^+$. We note $G(e_j)^-$ and $G(e_j)^+$ their
respective images in the vertex groups and $\vf_j:G(e_j)^-\longrightarrow G(e_j)^+$ the
isomorphism given by $\vf_j=\vf_j^+\circ (\vf_j^-)^{-1}$.
Now for an arbitrary edge $e$, $G(\bar{e})=G(e)$, $G(\bar{e})^-=G(e)^+$,
$G(\bar{e})^+=G(e)^-$, $\vf_{\bar{e}}^-=\vf_e^+$, $\vf_{\bar{e}}^+=\vf_e^-$, and
$\vf_{\bar{e}}=\vf_e^{-1}$.

Once a base point $*\in X$ is given, one defines the fundamental group $\pi_1(\cal{G},X,*)$ of
the pointed graph of groups $(\cal{G},X,*)$ ; it turns out that this definition does not depend
on $*\in X$ (cf. \cite{serre}) and that $\pi_1(\cal{G},X)$ is naturally isomorphic to $\P(M)$
(cf. \cite{scott-wall}). The vertex groups naturally embed in $\P(\cal{G},X)$ while edge groups
embed in vertex groups ; their respective images will be called {\sl vertex subgroups} of $\P(M)$
and {\it edge subgroups} of vertex groups (hence of $\P(M)$).

The splitting of $\pi_1(M)$ as a graph of groups provides a (finite) presentation, once a
maximal tree $T$ in $X$ is chosen. An edge of $X$ will be said to be $T${\it -separating}
if it belongs to $T$, and $T${\it -non separating} otherwise.
 If for each $i=1,2,\ldots ,n$ the vertex group $\P(M_i)$ admits the
(finite) presentation $<S_i\,|\, R_i>$ then $\P(M)$ admits the following presentation : \smallskip\\
{\bf Generators :} $\qquad S_1\cup S_2\cup\cdots\cup S_n\cup \{t_e\,|\, \text{$e$ is an edge
of $X$}\}$.\smallskip\\
{\bf Relators :}\begin{align*}\qquad R_1\cup R_2\cup \cdots \cup
R_n\; &\cup\; \{\text{for all edge $e$ of $X$}, \text{for all
$c\in G(e)^-$}, \quad t_e\, \vf_e(c)\,
t_e^{-1}=c\}\\
&\cup\; \{\text{for all edge $e$ of $X$},\;\;
t_{\bar{e}}=t_e^{-1}\}\\
 &\cup\; \{\text{for all $T$-separating
edge $e$},\;\; t_e=1\}
\end{align*}

A generator $t_e$
for some edge $e$, will be  called a {\it stable
letter} associated to the edge $e$ ; for an arbitrary edge $e_i$, we may also note $t_i$ instead of $t_{e_i}$.
Remark that since the edge groups are finitely generated, one
immediately obtains a finite presentation for $\P(M)$ by replacing
the relations $\forall c\in G(e)^-,\; t_e\, c\,
t_e^{-1}=\vf_e(c)$, by the same relations involving two
generators $c_1,c_2$ of $G(e)^-$.

\subsection{Algorithmically splitting the group}
Suppose the closed Haken manifold $M$ is given, in some manner,
for example by a triangulation, a Heegard splitting, or a Dehn
filling on a link. There exists an algorithm (cf. \cite{jt}),
which provides a JSJ decomposition of $M$ (as well as the
associated gluing maps). This algorithm uses an improvement of the
theory of normal surfaces due to Haken. It seeks until it finds a
maximal system of essential tori, which provides a JSJ
decomposition of $M$. Moreover this algorithm finds Seifert
invariants associated to each Seifert piece.

Once the JSJ decomposition, the pieces and the gluing maps are
given, one can easily split $\P(M)$ as a group of a graph of
groups, as described above. When we will be algorithmically
working on $\P(M)$, we will suppose a canonical presentation of
$\P(M)$ is given, that is a presentation of $\P(M)$ as above, such
that the Seifert pieces are given with their canonical
presentation (in such a way that they can be identified as being
Seifert pieces, and implicitly provide Seifert fibrations, cf. \S
4.1).

\section{\bf The conjugacy theorem}
As seen before the JSJ decomposition of a Haken closed manifold
$M$ provides a splitting of $\P(M)$ as a fundamental
group of a graph of groups $\P(M)=\P(\cal{G},X)$. This fact will
allow to establish a conjugacy theorem (i.e. which characterizes
conjugate elements) in $\P(M)$ in the spirit of analogous results
in amalgams or HNN extensions (cf. \cite{mks}, \cite{ls}). This is
the aim of this section. We will first need to recall classical
ways to write down an element of a group of a graph of groups, in
a "reduced form" before stating the main result, that is theorem
\ref{conjthm}.

\subsection{Cyclically reduced form} We first need to give some
definitions.

We note a path in $X$ in the extended manner :
$(v_{\s_1}, e_{\tau_1}, v_{\s_2}, e_{\tau_2}, \ldots
,v_{\s_m}, e_{\tau_m},v_{\s_{m+1}})$, where the $v_{\s_i}$
 and the $e_{\tau_i}$ are respectively vertices and edges of $X$ such that
$\ex(e_{\tau_i})=v_{\s_{i+1}}$ and $\o(e_{\tau_i})=v_{\s_i}$ for
$i=1,\ldots ,m$. The vertices $v_{\s_1}$ and $v_{\s_{m+1}}$ are its two endpoints ;
a loop is a path whom two endpoints co\"\i ncide.

Given an arbitrary graph of group $(\cal{G},X)$, we recall that in
the Bass-Serre's terminology (\cite{serre}) a {\sl word of type }
$\cal{C}$  is a
couple $(\cal{C},\mu)$ where :\\
-- $\cal{C}$ is a based loop in $X$, say
$\cal{C}=(v_{\s_1}, e_{\tau_1}, v_{\s_2}, e_{\tau_2}, \ldots
,v_{\s_m}, e_{\tau_m},v_{\s_{m+1}})$.\\
 -- $\mu$ is a sequence $\mu=(\mu_1,\mu_2,\ldots ,\mu_{m+1})$, such
that $\forall i=1,\ldots , m+1$,
 $\mu_i\in G(v_{\s_i})$ ; $\mu_i$ will be
called the {\it label} of the vertex $v_{\s_i}$.

The {\sl length } of a word of type $\cal{C}$ is defined to be the
length of the loop $\cal{C}$.
Once a base point $*$ in $X$ is given, each word of type $\cal{C}$
for some loop $\cal{C}$ with base point $*$, defines an element of
$\pi_1(\cal{G},X,*)$ (cf.\cite{serre}) which we will call its {\sl
label} ; the label of $(\cal{C},\mu)$ will be noted
$|\cal{C},\mu|$ and we shall speak of the {\sl form} $(\cal{C},\mu)$ for
$|\cal{C},\mu|$. When one considers a presentation as in \S 2.2,
$|\cal{C},\mu|=\mu_1.t_{\tau_1}.\mu_2.t_{\tau_2}\cdots
\mu_{m}.t_{\tau_m}.\mu_{m+1}$.

A word of type $\cal{C}$, $(\cal{C},\mu)$ is said to be {\sl
reduced} if either its length is 0 and its label is $\not= 1$,
or its length is greater than 1 and each time
$e_{\tau_{i-1}}=\bar{e}_{\tau_{i}}$ then $\mu_i\in
G(v_{\tau_i})\setminus G(e_{\tau_i})^+$. We shall speak of a {\sl
reduced form} for its label $|\cal{C},\mu|$.
Now given any non trivial element
$g\in \pi_1(\cal{G},X)$ there exists a reduced form associated
with $g$ (cf. \cite{serre}).\\

\noindent{\bf Definitions : } A {\sl cyclic conjugate} of
$(\cal{C},\mu)=((v_{\s_1},e_{\tau_1},v_{\s_2},\ldots,
e_{\tau_n},v_{\s_{n+1}}),(\mu_1,\mu_2,\ldots \mu_n,\mu_{n+1}))$
\nolinebreak is\linebreak[4]
$(\cal{C}',\mu')=((v_{\s_i},e_{\tau_i},v_{\s_{i+1}},\ldots,
e_{\tau_n},v_{\s_{n+1}}, e_{\tau_1},\ldots,
e_{\tau_{i-1}},v_{\s_i}),(\mu_i,\mu_{i+1},\ldots
\mu_n,\mu_{n+1}\mu_1,\ldots, \mu_{i-1},1))$, for some $i$, $0\leq
i\leq n$,
(indices are taken modulo $n$).\\
-- A word of type $\cal{C}$ is a {\sl cyclically reduced form}, if all of
its cyclic conjugates are reduced, and if $\mu_{n+1}=1$ when $n>1$ (hence this
last vertex becomes superfluous and should be forgotten).

One can associate to any non trivial conjugacy class in $\pi_1(\cal{G},X)$, a cyclically
reduced form whose label is an element of the class.
Just start from an element $g\not=1$ and find a reduced form for $g$ ; if it is not
cyclically reduced then reduce one of its cyclic conjugate. Apply this process as long as
it is possible ; since the length of the form strictly decreases it must stop providing a
cyclically reduced form whose label is a conjugate of $g$.

The crucial property for reduced forms is that given two reduced forms $(\cal{C},\mu)$ and
$(\cal{C},\mu')$ such that $|\cal{C},\mu|=|\cal{C},\mu'|$, then
necessarily $\cal{C}=\cal{C}'=(v_{\s_1},e_{\tau_1},v_{\s_2},\ldots ,e_{\tau_n},v_{\s_{n+1}})$
and there exists a sequence $(c_1,c_2,\ldots ,c_n)$, $c_i\in G(e_{\tau_i})$ such that
$\mu_1=\mu_1'.(c_1^-)^{-1}$, $\mu_{n+1}=c_n^+.\mu'_{n+1}$ in $G(v_{\s_{1}})$ and
$\mu_i=c_{i-1}^+.\mu_i'.(c_i^-)^{-1}$ in $G(v_{\s_i})$ for $i=2, 3\ldots,n$ (cf. \S 5.2, \cite{serre}).
The conjugacy theorem in the next section gives an analogous of this property when one
considers cyclically reduced forms and conjugacy classes instead of reduced forms and elements.

\subsection{The conjugacy theorem}

In this section, we will consider a graph of groups $(\cal{G},X)$ associated to
the JSJ decomposition of an arbitrary Haken closed manifold $M$.
We prove in this case the conjugacy theorem in the group of such
a graph of groups. Nevertheless the theorem remains true for any graph of groups : we prove
this result using a different method and in full generality in a work in preparation.

\begin{thm}
\label{conjthm}
 Suppose $(\cal{C},\mu)$ and $(\cal{C}',\mu')$ are
cyclically reduced forms, whose labels $\w$ and $\w'$ are
conjugate in $\P(M)$. Then, $(\cal{C},\mu)$ and $(\cal{C}',\mu')$
have the same length, and moreover, either :
\begin{itemize}
\item[$(i)$]
Their length is equal to $0$, $\cal{C}=\cal{C}'=(v_{\s_1})$, and
$\w$,
$\w'$ are conjugate in $G(v_{\s_1})$.\\
\item[$(ii)$] Their length is equal to $0$, and there exists a path
$(v_{\a_0},e_{\b_1},\ldots ,e_{\b_{p}},v_{\a_p})$ in $X$, and a
sequence $(c_1,c_2,\ldots , c_{p})$ with $\forall i=1,2,\ldots
,p$, $c_i$ lying in the edge group $G(e_{\b_i})$, such that $\w\in
G(v_{\a_0})$, $\w'\in G(v_{\a_p})$, and
\begin{gather*}
\w\sim c_1^-\quad \text{in $G(v_{\a_0})$}\\
\w'\sim c_{p}^+\quad \text{in $G(v_{\a_p})$}\\
\text{and \;$\forall\, i=1,2,\ldots ,p-1$,}\quad c_i^+\sim
c_{i+1}^-\quad \text{in $G(v_{\a_i})$}
\end{gather*}\smallskip
\item[$(iii)$] Their length is greater than $0$. Up to cyclic
permutation of $(\cal{C}',\mu')$, the loops $\cal{C}$, $\cal{C}'$ are equal,
$\cal{C}=\cal{C}'=(v_{\s_1},e_{\tau_1},\ldots
,v_{\s_n},e_{\tau_n})$, and there exists a sequence $(c_1,\ldots
,c_n)$, with for all $i=1,\ldots ,n$, $c_i$ lying in the edge
group $G(e_{\tau_i})$, such that :
\begin{gather*}
\mu_1=c_n^+.\mu_1'.(c_1^-)^{-1}\quad \text{in $G(v_{\s_1})$}\\
\forall\, i=2,3,\ldots , n \quad
\mu_i=c_{i-1}^+.\mu_i'.(c_i^-)^{-1}\quad \text{in $G(v_{\s_i})$}
\end{gather*}
in particular, the element $c_n^+\in G(e_{\tau_n})^+$ conjugates
$\w'$ into $\w$ in $\P(M)$ :
$$\w=c_n^+.\w'.(c_n^+)^{-1}\quad \text{in $\P(M)$}$$
\end{itemize}
(Recall that if $c$ lies in the edge group $G(e_i)$, we note
$c^-=\vf_i^-(c)\in G(e_i)^-$ and $c^+=\vf_i^+(c)\in G(e_i)^+$.)
\end{thm}

\noindent{\bf Proof.}
Consider two cyclically reduced forms
$(\cal{C},\mu)$ and $(\cal{C},\mu')$ with respective labels $\w$ and $\w'$.

In order to define edge groups, vertex subgroups and the embeddings  we have considered
in \S 2.2 one base point $*_i$ in each piece $M_i$, one base point $\ul{*}_j$ in each torus
component of $W$, and for each edge $e_j$ with $\o(e_j)=v_k$, $\ex(e_j)=v_l$ two paths that
we note $[*_k,\ul{*}_j^-]$ and  $[*_l,\ul{*}_j^+]$ respectively in $M_k$ from $*_k$ to $\ul{*}_j^-$
and in $M_l$ from $*_l$ to $\ul{*}_j^+$. We deform by homotopy keeping their endpoints fixed
all of these paths in such a way that once we have noted  $[e_j]=[*_k,\ul{*}_j^-].[*_l,\ul{*}_j^+]^{-1}$,
all such $[e_j]$ become smooth and transverse with $W$ in $M$.

One constructs a smooth based
loop in $M$ representing $(\cal{C},\mu)$ in the following way : suppose
$\cal{C}=(v_{\s_1}, e_{\tau_1}, v_{\s_2}, e_{\tau_2}, \ldots ,
e_{\tau_n})$ ; for each vertex $v_{\s_i}$ of $\cal{C}$ choose a
smooth loop $V_i$ with base point $*_{\s_i}$ in $int (M_{\s_i})\subset
M$, which represents the label $\mu_i$ of $v_{\s_i}$ in
$\P(M_{\s_i},*_{\s_i})$. Replace in $\cal{C}$ the vertex $v_{\s_i}$ by this
loop. Replace each edge $e_{\tau_j}$ in $\cal{C}$ by the path $\lb
e_{\tau_j}\rb$ of $M$. Finally, concatenate the elements of the
sequence obtained and deform by small $*_{\s_i}$-homotopy each $V_i$ to obtain a
smooth path $\cal{P}_\w$, with
base point $*_{\s_1}$.
Proceed in the same way to obtain a smooth based
path $\cal{P}_{\w'}$ representing $(\cal{C}',\mu')$.\\

Suppose that $\w$ and $\w'$ are conjugate in $\P(M)$. Then the
loops $\cal{P}_\w$ and $\cal{P}_{\w'}$ are freely homotopic in
$M$. Hence, there exists a map $f:S^1\times I\longrightarrow M$,
such that $f$ restricted to $S^1\times 0$ is $\cal{P}_\w$, and $f$
restricted to $S^1\times 1$ is $\cal{P}_{\w'}$. One can also
suppose that $f$ is smooth. Since $\cal{P}_\w$ and $\cal{P}_{\w'}$
are transverse to $W$, according to the homotopy transversality
theorem one can deform $f$ without changing either $\cal{P}_w$ or
$\cal{P}_{\w'}$, such that it becomes transverse to $W$. Then, by
the transversality theorem, $f^{-1}(W)$ is a compact 1-submanifold
of $S^{1}\times I$, such that
$\partial(f^{-1}(W))=f^{-1}(W)\cap\partial(S^1\times I)$ ; hence $f^{-1}(W)$
consists of disjoint segments and circles properly
embedded in $S^1\times I$.
Among all the ways to choose and deform $f$ as above, we consider one such that the map $f$
obtained after deforming is { minimal}, in the sense that the number of connected components
of $f^{-1}(W)$ is minimal.

The minimality in the choice of $f$ implies that none of the
circle components of $f^{-1}(W)$ bound a disk, while the facts
that $(\cal{C},\mu)$ and $(\cal{C}',\mu')$ are cyclically reduced
forms implies that none of the segment components of $f^{-1}(W)$
have its two boundary components both in $S^1\times 0$ or in
$S^1\times 1$. Hence if $f^{-1}(W)$ is non-empty, then it consists
either of disjoint circles parallel to the boundary or of disjoint
segments joining $S^1\times 0$ to $S^1\times 1$ (cf. figure 2).
\begin{figure}[ht]
\center{\includegraphics[scale=0.5]{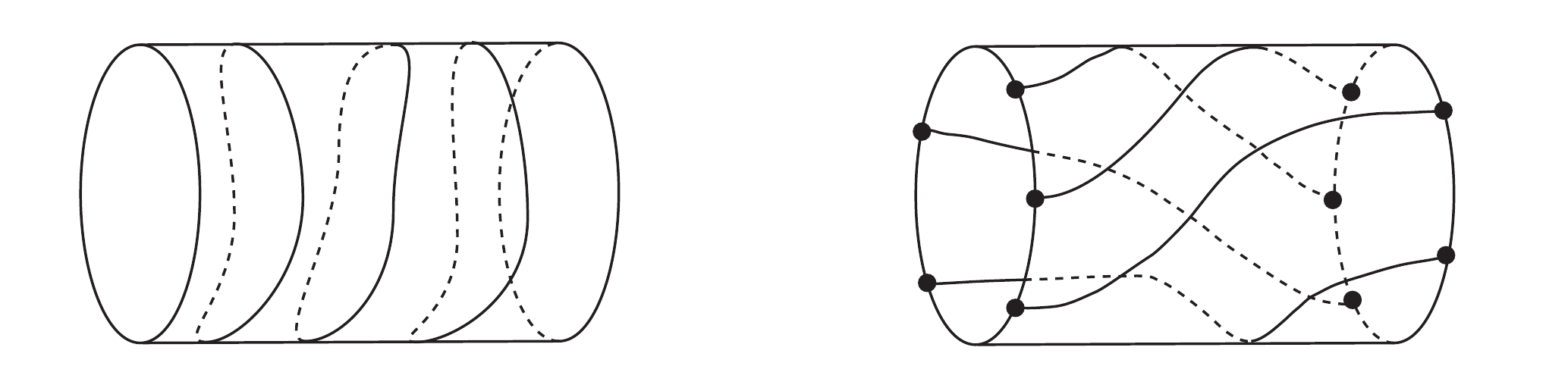}} \caption{}
\end{figure}

First we can conclude that $\cal{P}_\w$ and $\cal{P}_{\w'}$
intersect $W$ the same number of times. So the cyclically reduced forms
$(\cal{C},\mu)$ and $(\cal{C}',\mu')$ must have the same length.\\

\noindent{\bf First case :} Suppose that $f^{-1}(W)$ is empty. So
the annulus $f(S^1\times I)$ lies in some piece, say $int(M_i)$.
It implies that $(\cal{C},\mu)$, $(\cal{C}',\mu')$ both have length 0, and moreover that
$\w,\w'$ belong to $\P(M_i)$ and are conjugate in $\P(M_i)$ ; hence conclusion $(i)$ holds.\\

\noindent {\bf Second case :} Suppose $f^{-1}(W)$ consists of $p$
circles, $C_1,C_2,\ldots ,C_p$. In this case $(\cal{C},\mu)$ and $(\cal{C}',\mu')$ both have
length equal to 0. We note $C_0=S^1\times 0$, $C_{p+1}=S^1\times 1$, and consider
the circles $C_i$ as based loops such that they are all isotopic in $S^1\times I$ and
$f\circ C_0=\cal{P}_\w$, $f\circ C_{p+1}=\cal{P}_{\w'}$. We will proceed by induction on $p$
to show that conclusion $(ii)$ holds.

Consider first the case $p=1$. If one cuts $S^1\times I$ along $C_1$ it decomposes in two annuli
bounded respectively by $C_0,C_1$ and $C_1,C_2$ which map respectively under $f$, say in $M_k$,
$M_l$ ; note $\T_i$ the component of $W$ in which $C_i$ maps and $e_i=e(\T_i)^{\pm 1}$ the edge
such that $\o(e_i)=v_k$, $\ex(e_i)=v_l$. The loop $f\circ C_0=\cal{P}_\w$ (resp. $f\circ C_{2}=\cal{P}_{\w'}$)
represents the element $\w$ in $G(v_k)$ (resp. $\w'$ in $G(v_l)$). The loop $f\circ C_1\subset \T_i$
defines a conjugacy class $[c_1]$ in $G(e_i)$ such that $\w\sim c_1^-$ in $G(v_k)$ and
$\w'\sim c_1^+$ in $G(v_l)$. Hence conclusion $(ii)$ holds when one considers the path $(v_k,e_i,v_l)$
and the sequence $(c_1)$.

\begin{figure}[ht]
\center{\includegraphics{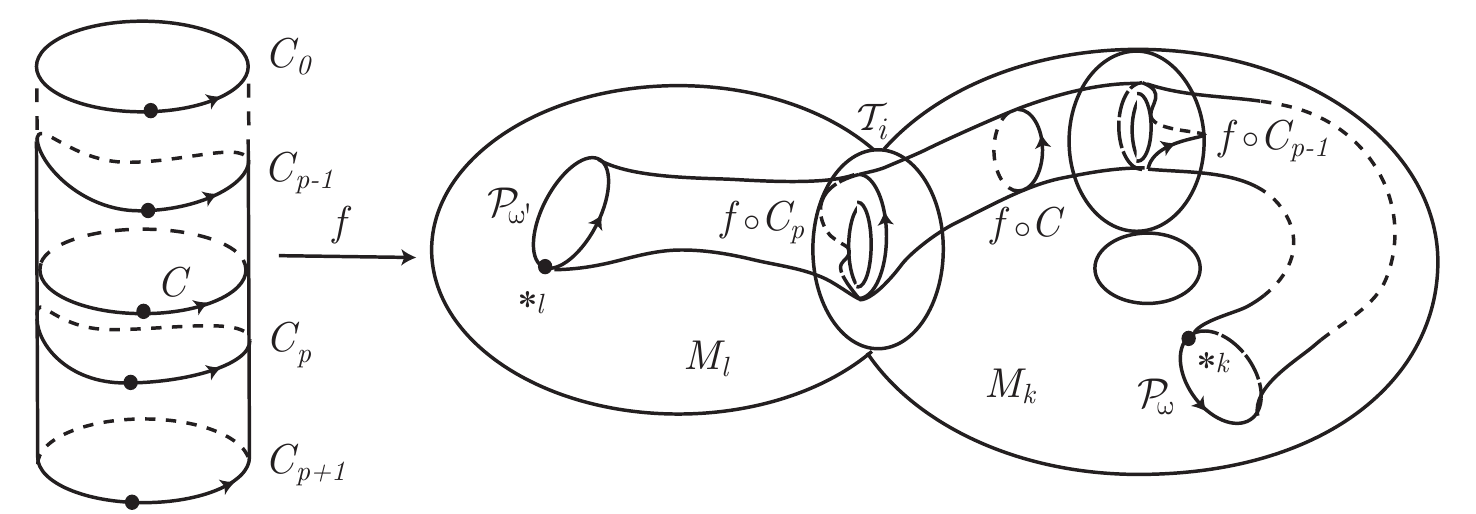}} \caption{}
\end{figure}

Consider now the case $p>1$. Suppose that conclusion $(ii)$ holds whenever $f^{-1}(W)$ consists of $p-1$
circles, and moreover that $f^{-1}(W)$ has $p$ components. The loops $C_{p-1}$ and $C_p$ cobound an annulus
$A$ which maps in, say $M_k$ ; note $\T_i$ the component of $W$ in which $C_p$ maps. Consider an additional
loop $C$ in $int(A)$ isotopic in $A$ with both $C_{p-1}$ and $C_p$. Then $f\circ C$ defines a conjugacy
class $[\w_0]$ in $G(v_k)$, and once $e_i=e(\T_i)^{\pm 1}$ is judiciously chosen, $f\circ C_p$ defines a
conjugacy class $[c_p]$ in $G(e_i)$ such that $c_p^-\sim \w_0$ in $G(v_k)$.  Then cut $S^1\times I$ along
$C$ : it decomposes in two parts, the former one $A_0$ containing $S^1\times 0$ and the latter $A_1$
containing $S^1\times 1$. The hypothesis of induction can be applied when one restricts $f$ to the annulus
$A_0$ to provide a path, say $\cal{P}=(v_0,e_1,\ldots ,e_{p-1},v_k)$ with endpoint $v_k$ and a sequence
$c=(c_1,\ldots ,c_{p-1})$ as in conclusion $(ii)$ from $\w$ to $\w_0$. The same argument as in the former
case $p=1$ applied to the annulus $A_1$ provides the path $(v_k,e_i,\ex(e_i))$ and the sequence $(c_p)$
from $\w_0$ to $\w'$. Then conclusion $(ii)$ holds when one considers the path $\cal{P}.(v_k,e_i,\ex(e_i))$
together with the sequence $c.(c_p)$.\\

\noindent {\bf Third case :} Suppose $f^{-1}(W)$ consists of $n>0$
segments ; then $n$ is the length of both $(\cal{C},\mu)$ and
$(\cal{C}',\mu')$. We must show that conclusion $(iii)$ holds. We note $\cal{C}=(v_{\s_1},
e_{\tau_1}, v_{\s_2}, e_{\tau_2}, \ldots , e_{\tau_n})$, while $\cal{C}'$ is clearly a cyclic
conjugate of $\cal{C}$.
For more convenience during the rest of the proof indices  will be given modulo $n$.
We note $x_1,x_2,\ldots ,x_n$ the points of
$\cal{P}_\w^{-1}(W)$ in such a way that if one starts from
$(1,0)$ and turns in the positive sense on $S^1\times 0$,
one  meets $x_1,x_2,\ldots x_n$ in this order, and we proceed the
same way with the points $x_1',x_2',\ldots ,x_n'$ of
$\cal{P}_{\w'}^{-1}(W)$ ; they decompose $S^1\times 0$ and $S^1\times 1$ in  paths which will
be respectively noted  $[x_{i-1},x_{i}]$ and $[x_{i-1}',x_{i}']$, $i=1,2,\ldots n$.
Consider
the segments of $f^{-1}(W)$ as paths $C_1,C_2,\ldots ,C_n$ such
that each $C_i$ starts in $x_i\in S^1\times 0$ and ends in some $x_j'\in S^1\times 1$.
Necessarily there exists an integer $p$ such that $\forall\,
i=1,\ldots ,n$, the path $C_i$ ends in $x_j'\in S^1\times 1$ with
$j=i+p$. By changing if necessary $(\cal{C}',\mu')$ into a cyclic conjugate, we can suppose that $p=0$ ;
hence with this convention $\cal{C}=\cal{C}'$ and the paths $C_i$ go from $x_i$ to $x_i'$. The annulus
decomposes into $n$ strips such that the boundary of strip $i$ contains the loop
$[x_{i-1},x_{i}].C_{i}.[x_{i-1}',x_{i}']^{-1}.C_{i-1}^{-1}$.

By construction each of the $C_i$ maps under $f$ on a loop in $\T_{\tau_i}$ with base point
$f(x_i)=f(x_i')=\ul{*}_{\tau_i}$. Hence, $f\circ
C_i$ defines an element $c_i\in G(e_{\tau_i})$.
Consider for $i=1,2,\ldots ,n$, the loops ${C}_i^-$ and
${C}_i^+$ defined by :
\begin{gather*}
{C}_i^-=[*_{\s_i},\ul{*}_{\tau_i}^-].f\circ{C}_i.[\ul{*}_{\tau_i}^-,*_{\s_i}]\\
{C}_i^+=[*_{\s_{i+1}},\ul{*}_{\tau_i}^+].f\circ{C}_i.[\ul{*}_{\tau_i}^+,*_{\s_{i+1}}]
\end{gather*}
The loop ${C}_i^-$ (resp. ${C}_i^+$) has base point $*_{\s_i}$ (resp. $*_{\s_{i+1}}$) and represents the
element $c_i^-\in G_{\tau_i}^-\subset\P(M_{\s_i})$ (resp. $c_i^+\in G_{\tau_i}^+\subset\P(M_{\s_{i+1}})$),
with $c_i^+=\vf_i^+(c_i)$ et $c_i^-=\vf_i^-(c_i)$.

Consider also for $i=1,2,\ldots, n$ the loops $W_i$ and $W_i'$ in
$int(M_{\s_i})$, with base point $*_{\s_i}$ defined by :
\begin{gather*}
W_i=[*_{\s_i},\ul{*}_{\tau_{i-1}}^+].f\circ[x_{i-1},x_i].[\ul{*}^-_{\tau_i},*_{\s_i}]\\
W_i'=[*_{\s_i},\ul{*}_{\tau_{i-1}}^+].f\circ[x_{i-1}',x_i'].[\ul{*}^-_{\tau_i},*_{\s_i}]
\end{gather*}
By construction, once we have noted $\mu=(\mu_1,\mu_2,\ldots,\mu_n)$ and $\mu'=(\mu_1',\mu_2',\ldots,\mu_n')$,
the loops
$W_i$ and $W_i'$  represent respectively the elements $\mu_i$ and
$\mu_i'$ of $\P(M_{\s_i})$, .

\begin{figure}[ht]
\center{\includegraphics{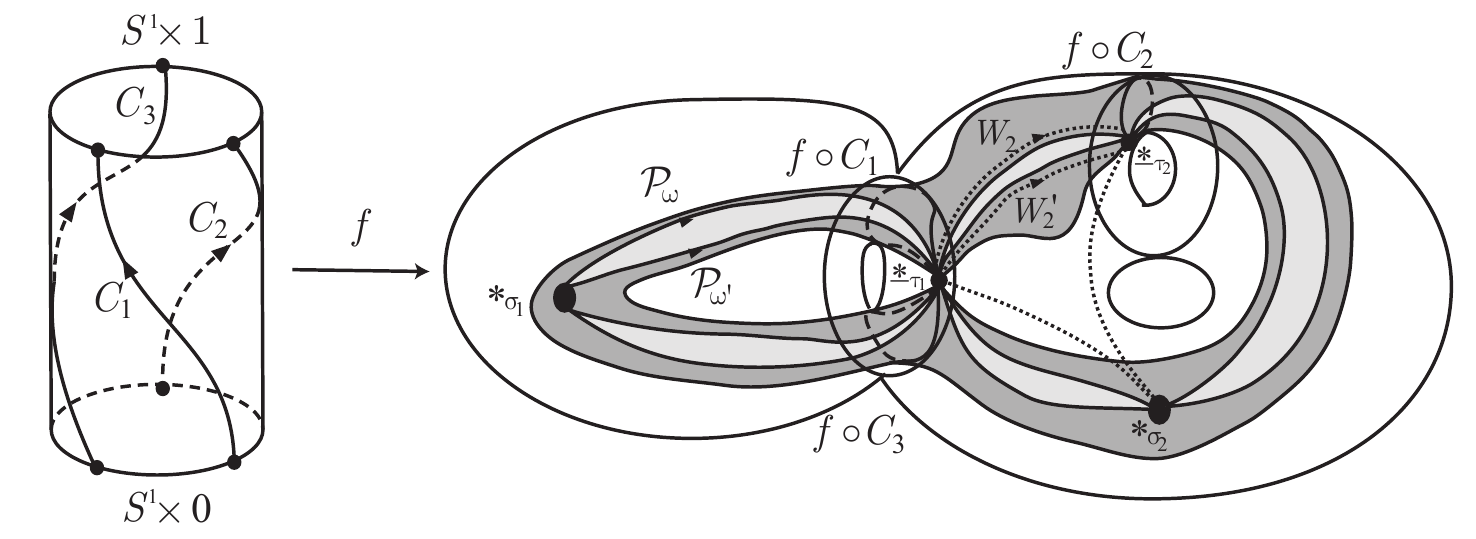}} \caption{}
\end{figure}

The strips show that for $i=1,2,\ldots,n$ the paths $f\circ [x_{i-1},x_i]$
and $f\circ C_{i-1}.f\circ [x'_{i-1},x'_i].f\circ C_{i}^{-1}$ are homotopic in $M_{\s_i}$
 with endpoints $\ul{*}_{\tau_{i-1}}$, $\ul{*}_{\tau_i}$ fixed. Hence  for all $i=1,2,\ldots ,n$,
 one has in $M_{\s_i}$ the $*_{\s_i}$-homotopy :
\begin{align*}
{C}_{i-1}^+.W_i'.({C}_i^-)^{-1} &\approx_{*_{\s_i}}
[*_{\s_i},\ul{*}_{\tau_{i-1}}^+].f\circ{C}_{i-1}.f\circ[x_{i-1}',x_i'].f\circ{C}_i^{-1}.[\ul{*}_{\tau_i}^-,*_{\s_i}]\\
 &\approx_{*_{\s_i}}
[*_{\s_i},\ul{*}_{\tau_{i-1}}^+].f\circ[x_{i-1},x_i].[\ul{*}_{\tau_i}^-,*_{\s_i}]\\
&\approx_{*_{\s_i}} W_i
\end{align*}
and one can slightly deform the paths on  regular
neighborhoods of $\cal{T}_{\tau_{i-1}}$ and $\cal{T}_{\tau_i}$
such that the homotopy takes place in $int(M_{\s_i})$.  Hence, for all
$i=1,2,\ldots ,n$,
$\mu_i=c_{i-1}^+.\mu_i'.(c_i^-)^{-1}$ in $\P(M_{\s_i})$.

Now,
\begin{align*}
\w &=\mu_1.t_{\tau_1}.\mu_2.t_{\tau_2}.\cdots t_{\tau_{n-1}}.\mu_n.t_{\tau_n}\\
\w &=c_n^+.\mu_1'.(c_1^-)^{-1}.t_{\tau_1}.c_1^+.\mu_2'.(c_2^-)^{-1}.t_{\tau_2}.c_2^+\cdots (c_{n-1}^-)^{-1}.t_{\tau_{n-1}}.c_{n-1}^+.\mu_n'.(c_n^-)^{-1}.t_{\tau_n}\\
\w &= c_n^+.\mu_1'.t_{\tau_1}.\mu_2'.t_{\tau_2}.\cdots t_{\tau_{n-1}}.\mu_n'.t_{\tau_n}.(c_n^+)^{-1}\\
\end{align*}
shows that
$\w=c_n^+.\w'.(c_n^+)^{-1}$ in $\P(M)$, with $c_n^+\in
G(e_{\tau_n})^+$. So in this case the conclusion $(iii)$ holds,
which concludes the proof.\hfill\bs

\vskip 0.4cm

\section{\bf Reducing conjugacy problem in $\P(M)$ to problems in the pieces}

We establish in this section the main argument for solving the
conjugacy problem in the group of a Haken closed manifold $M$
(which is not a $S^1\times S^1$-bundle or two twisted $I$-bundles
over $\KB_2$ glued along their boundary). We reduce the conjugacy
problem to three elementary problems in the group of the pieces
obtained in a non trivial JSJ decomposition of $M$ : namely the conjugacy problem (of course), the boundary
parallelism problem and the 2-cosets problem.

The {\it boundary parallelism problem} consists in, given a
boundary subgroup $T$ of the group $\P(N)$ of a piece $N$, to
decide for any element $\w$ of $\P(N)$ (given as words on a given
set of generators), if $\w$ is conjugate in $\P(N)$ to an element
of $T$.

The {\it 2-cosets problem} consists in, given two boundary
subgroups $N_1,N_2$ of $\P(N)$ (possibly identical),  to find for
any $u,v\in \P(N)$, all the couples $(c_1,c_2)\in N_1\times N_2$
which are  solutions of the equation $u=c_1.v.c_2$ in $\P(N)$.

We show that if one can solve those three problems in the groups
of the pieces obtained, then one can solve the conjugacy problem
in the group of $M$ (theorem \ref{reduce}).
Remark first that we
will suppose that a canonical presentation of $\P(M)$, that is its
decomposition as a graph of groups, as well as canonical
presentations for the groups of the Seifert pieces, are given.
Indeed, given the manifold $M$, there exists an algorithm based
along the same lines as the Haken theory of normal surfaces, which
provides a minimal JSJ decomposition of $M$, as well as fibrations
of the Seifert pieces (\cite{jt}). Moreover, given a finite
presentation of the group of a Haken manifold $M$, one can
 reconstruct a triangulation of the manifold $M$.

In the previous section we established a theorem characterizing
conjugate elements, which does not directly provide a  solution to
the conjugacy problem, but which is essential to the reduction. On
its own, this result does not allow such a reduction, but the key
point is that groups of the pieces of a JSJ decomposition have
algebraic properties (a kind of "malnormality" for the boundary
subgroups, proposition \ref{annuli}) which together with the lemma
\ref{fibers} (a consequence of the minimality of the JSJ
decomposition), make the reduction process work. These algebraic
properties will be established in \S 4.2, and will imply the
$k$-acylindricity of the JSJ splitting (\S 4.3), as well as the
existence of an algorithm to write words in cyclically reduced
forms (\S 4.4), which are all essential to the reduction process
(\S 4.5). But first, we  recall some elementary facts upon Seifert
fibre spaces (we refer the reader to \cite{seifert}, \cite{js},
\cite {jaco}, \cite{orlik}).

\subsection{Reviews on Seifert fiber spaces}
Let $M$ be a Seifert fibered space. A Seifert fibration of $M$ is
characterised by a set of invariants (up to fiber preserving
homeomorphism) of one of the forms :
\begin{gather*}
(o,g,p,b\mid \a_1,\b_1,\a_2,\b_2,\ldots ,\a_q,\b_q)\\
(n,g,p,b\mid \a_1,\b_1,\a_2,\b_2,\ldots ,\a_q,\b_q)
\end{gather*}
The former case occurs when the base is oriented ("$o$" stands for
"oriented"), and the latter when the base is non-oriented ("$n$"
for "non-oriented"). The numbers "$g$", "$p$" are respectively the
genus of the base, and the number of its boundary components. The
number $b$ is related to the Euler number of the $S^1$-bundle
associated to the fibration, and  $q$ is the number of exceptional
fibers ; $\a_i$ is the index of the $i$-st exceptional fiber,
$0<\b_i<\a_i$ and $(\a_i,\b_i)$ is the type of this exceptional
fiber.

A Seifert fiber space may admit several fibrations, but it remains
 isolated cases : at the exception of lens spaces, prism
 manifolds, a solid torus, a twisted $I$-bundle over $\KB_2$, or
 the double of a twisted $I$-bundle over $\KB_2$, a Seifert fiber
 spaces can be endowed with a unique Seifert fibration (\cite{jaco}, theorem VI-17).

Now, given a set of invariants of its Seifert fibration, $\P(M)$
admits a canonical presentation, of one of the forms, according to
whether its base is oriented or non-oriented (cf. \cite{jaco}).
\begin{align*}
<a_{1},& b_{1},\ldots,a_{g},b_{g},c_{1},\ldots,c_{q}, d_{1}
,\ldots,d_{p},h \mid\\
& \lbrack a_{i},h\rbrack =\lbrack b_{i},h\rbrack =\lbrack
c_{j},h\rbrack=\lbrack d_{k},h\rbrack=1; c_{j}^{\alpha
_{j}}=h^{\beta _{j}}; h^{b}=(\prod_{i=1}^{g} \lbrack
a_{i},b_{i}\rbrack)c_{1}
\cdots c_{q}d_{1}\cdots d_{p}>\qquad (1)\\
<a_{1},& \ldots ,a_{g},c_{1},\ldots,c_{q},d_{1},\ldots, d_{p},h\mid\\
& a_{i}ha_{i}^{-1}=h^{-1} ;\lbrack c_{j},h\rbrack=\lbrack
d_{k},h\rbrack=1;c_{j}^{\alpha _{j}}=h^{\beta _{j}};
h^{b}=(\prod_{i=1}^g a_{i}^{2})c_{1}\cdots c_{q}d_{1}\cdots d_{p}>
\quad\quad\quad\quad\;\, (2)
\end{align*}
with $1\leq i\leq g$, $1\leq j\leq q$, and $1\leq k\leq p$. The
generator $h$ is the class of a (any) regular fiber, and if
$\cal{T}_k$ is a component of $\partial M$, the associated
boundary subgroup $T_k$ is generated by $h,d_k$. The element $h$
generates a normal subgroup $N=<h>$ of $\P(M)$, called the {\sl
fiber}. Moreover, if $\P(M)$ is infinite, then $h$ has infinite
order (lemma II.4.2, \cite{js}), and hence, there is  an exact
sequence :
$$1\longrightarrow \Z \longrightarrow \P(M) \longrightarrow \P(M)/N
\longrightarrow 1$$ By looking at the presentation above, we see
that if the base is oriented, $N$ is central.

When the base is non-oriented, let $C$ be the subgroup of $\P(M)$,
of all elements $\w$ written as words on the canonical generators with
an even number of occurrences of generators
$a_1,a_2,\ldots ,a_g$ and their inverses. According to the
relators of $\P(M)$, obviously, such a fact does not depend on the
word chosen in the class of $\w$. One easily shows that $C$ has
index 2 in $\P(M)$, and that $C$ is the centralizer of any
non-trivial element of $N$, while, for all $u\not\in C$,
$u.h.u^{-1}=h^{-1}$. When the base is oriented, just set $C=\P(M)$
; obviously, $C$ is the centralizer of any element of $N$. This
combinatorial  definition of the subgroup $C$, agrees with the
topological definition of the {\it canonical
subgroup} of $\P(M)$, as seen in \cite{js}.\\

Let us focus --as an example-- on the $I$-twisted bundle over
$\KB_2$ (that we shall call $K$) in order to recall elementary
facts needed later. The group of $K$ is the group of the Klein
bottle \mbox{$\P(K)=<a,b\,|\,a.b.a^{-1}=b^{-1}>$}. Its boundary
consists of one toro\" \i dal component, and the boundary subgroup
is the (free abelian of rank 2) subgroup of index 2 of $\P(K)$ :
$\P(\partial K)=<a^2,b>$.

One can endow $K$ with two Seifert fibrations. The first has base
the M\"obius band, and no exceptional fiber. The class of a
regular fiber is $b$. In this case $N$ is not central, and the
canonical subgroup is $<a^2,b>$. The second seifert fibration has
base a disk, and two exceptional fibers of index 2. The class of a
regular fiber is $a^2$, $N$ is central, and the canonical subgroup
is the whole group $\P(K)$.

\subsection{Algebraic properties in the pieces}
To proceed  we first need to establish an
 algebraic property  which is essential to  the reduction.
Recall that if $M$ is a manifold with non-empty boundary,
 and $\cal{T}$ is a connected component of $\partial M$, the canonical embedding
 $i :\cal{T}\hookrightarrow M$
 defines a conjugacy class of subgroups of $\P(M)$ : the subgroup $T=i_\ast(\P(\cal{T}))$
  depends of the  choice of a path from
the base point of $M$ to the base point of $\cal{T}$. Choosing
another such path changes $T$ in $gTg^{-1}$ for some $g\in \P(M)$.
Each element of the conjugacy class of $T$, will be called a {\it
boundary subgroup} of $\P(M)$ associated to $\cal{T}$.

Remark that one easily verifies that if $M$ is Haken and not
homeomorphic to a thickened surface, then boundary subgroups
associated to distinct boundary components are non conjugate.

\begin{prop}
\label{annuli}
Suppose $M$ is a piece obtained in a non trivial JSJ decomposition
of a Haken closed manifold
 which is not a $S^1\times
S^1$-bundle over $S^1$. Fix for each component of $\partial M$ a
boundary subgroup $T_i$.
\begin{itemize}
\item{} If $M$ is hyperbolic, and $T_1,T_2$ are two non conjugate
boundary subgroups of $\P(M)$, then no non trivial
element of $T_1$ is conjugate in $\P(M)$ to an element of $T_2$.
For any boundary subgroup $T_1$, if two elements $t,t'\in T_1$ are
conjugate by an element $u\in \P(M)$, then necessarily $t=t'$ and
$u\in T_1$.

\item{} If $M$ is a Seifert fibered space, and is not the twisted $I$-bundle over the Klein bottle,
it admits a unique fibration. We call $h$ the class in $\P(M)$ of
a regular fiber and $C$ the canonical subgroup. If $T_1,T_2$ are two non conjugate boundary
subgroups, then $<h>\subset T_1\cap T_2$, and if $v\in \P(M)$
conjugates $t_1\in T_1$ into $t_2\in T_2$, then $t_1,t_2\in <h>$
and $t_1=t_2^{\pm 1}$, with $t_1=t_2$ precisely when $v\in C$. For
any boundary subgroup $T_1$, if $t,t'\in T_1$ are conjugate by an
element $u\in \P(M)$, then either $t,t'\in <h>$ and $t'=t^{\pm 1}$
with $t=t'$ exactly when $u\in C$, or $t=t'$ and $u\in T_1$.
\item{} If $M$ is the twisted $I$-bundle over the Klein bottle,
then $\P(M)=<a,b|aba^{-1}=b^{-1}>$, and $\P(\partial M)=<a^2,b>$.
Two elements of $\P(\partial M)$, $a^{2n}b^p$ and $a^{2m}b^q$ are
conjugate in $\P(M)$, if and only if $n=m$ and $p=\pm q$.
\end{itemize}
\end{prop}

\noindent {\bf Proof.} Remark first that $M$ is a boundary
irreducible Haken manifold, with non-empty boundary consisting of
tori. So $M$ cannot be homeomorphic to $S^1\times D^2$. Moreover,
according to the lemma \ref{s1s1i},
 $M$ is not homeomorphic to $S^1\times S^1\times I$.

 The third case is easy to check from the presentation given. We leave it as an exercise
 for the reader. We will prove the two remaining cases
separately.

Suppose first that $M$ is hyperbolic. Its fundamental group
$\P(M)$ is a torsion free discrete subgroup of $PSL(2,\C)$, which
acts by isometries on $\H^3$. This action naturally extends to
$\H^3\cup\partial \H^3$. Each boundary subgroup $T_i$ corresponds
to a  maximal parabolic subgroup of $\P(M)$, with limit point the
cusp point $p_i\in\partial{H}^3$ (i.e. each element of $T_i$ is
parabolic and fixes $p_i$ and conversely each parabolic element
which fixes $p_i$ is in $T_i$ ; cf. \cite{rat}, \S 12.2). Suppose
that $u$ in $\P(M)$ conjugates two elements of $T_1$ ; then  $u$
must fix $p_i\in\partial\H^3$. According to the theorem 5.5.4 of
\cite{rat}, $u$ cannot be loxodromic, and hence is parabolic, so
$u\in T_1$, which proves the second part of the assertion.

Now suppose $T_1,T_2$ are two distinct boundary subgroups,
characterized by two (distinct) cusp points $p_1$ and $p_2$. If an
element $u$ conjugates two non trivial elements of $T_1$ and
$T_2$, then $u.p_1=p_2$. So, the fact that this cannot occur, is a
direct implication of the well-known fact that connected
components of $\partial M$ are in 1-1 correspondence with orbits
under the action of $\P(M)$ of the set of cusp points (cf. \cite{rat} \S 12.2),
which concludes the proof in this  case.\\

Suppose now that $M$ is a Seifert fibered space. By hypothesis,
$M$ has non-empty boundary, and is neither $D^2\times I$, nor the
twisted $I$-bundle over $\KB_2$, and hence admits a unique
fibration. Then $<h>$ is an infinite cyclic normal subgroup of
$\P(M)$, which does not depend on the regular fiber
 considered (cf. lemma
II.4.2, \cite{js}). Moreover, any component $\cal{T}$ of $\P(M)$
is trivially fibered by regular fibers. Hence any (free abelian of
rank two) boundary subgroup of $\P(M)$ contains $<h>$ as a
subgroup, which proves the beginning of the assertion.

Now suppose that two non trivial elements $t_1,t_2$ in the
respective  boundary subgroups $T_1,T_2$ (possibly $T_1=T_2$) are
conjugate in $\P(M)$. This gives rise to a map of pairs $f:
(S^1\times I,
\partial (S^1\times I)) \longrightarrow (M,\partial M)$ such that
its restrictions on $S^1\times 0$ and $S^1\times 1$ are non
contractible loops $\tau_1$, $\tau_2$ representing respectively
$t_1$ and $t_2$.

If this map is essential, according to the lemma II.2.8 of
\cite{js}, $\tau_1$, $\tau_2$ are homotopic in $\partial M$ to
powers of regular fibers, hence $t_1,t_2\in\ <h>$. Now, it appears
clearly from the presentation of $\P(M)$ (cf. \S 4.1), that for
all $u\in \P(M)$, $uhu^{-1}=h^{\e}$, with $\e=\pm 1$. Hence, in
this case $t_1=t_2^{\e}$ for some $\e=\pm 1$. Moreover, $C$ is the
centralizer of any non trivial element of $<h>$, and so with these
notations, $\e=1$ precisely when $u\in C$.

If the map is not essential, $f$ is homotopic rel.\ $\partial
(S^1\times I)$, with a map $g : S^1\times I\longrightarrow
\partial M$. Necessarily, $\tau_1$ and $\tau_2$ are in the same
component of $\partial M$, and homotopic in this component. Hence,
$T_1=T_2$ and $t_1,t_2$ are conjugate in $T_1$, but  since $T_1$
is abelian, $t_1=t_2$, which concludes the proof of the
assertion.\hfill \bs

\subsection{Acylindricity of the JSJ decomposition}
In this section, $M$ stands for a haken closed manifold which is
neither a $S^1\times S^1$-bundle over $S^1$, nor obtained by
gluing two twisted $I$-bundles over $\KB_2$.

 As a first
consequence of the lack of annuli stated in the previous
paragraph, we establish an essential result which asserts the
acylindricity (in the sense of Sela, \cite{acyl}) of the
Bass-Serre tree associated to the JSJ decomposition.  Roughly
speaking it means that there exists $K>0$, such that if two curves
lying in two pieces are freely homotopic in $M$, then any homotopy
between them can be deformed to meet the JSJ surface at most $K$
times.

Suppose $(\cal{G},X)$ is the JSJ splitting of $\P(M)$, and
$u,v\in\P(M)$ are conjugate elements lying in vertex groups.
According to the theorem \ref{conjthm} $(i)$, $(ii)$, there exists a
path $(v_{\a_0},e_{\b_1},\ldots,e_{\b_p},v_{\a_p})$ ($p\geq 0$), and a
sequence $(c_1,c_2,\ldots c_p)$ following the conclusion of the
theorem. We will say this sequence is {\it reduced} if whenever
$e_{\b_i}=\bar{e}_{\b_{i+1}}$, then $c_i^+$ and $c_{i+1}^-$ are non
conjugate in $G(e_{\b_i})^+$ (otherwise the sequence can be
shortened). The integer $p$ is called the {\it length} of the
sequence.

We will say that the JSJ splitting is \textit{$k$-acylindrical},
if whenever $u,v\in\P(M)$ are conjugate elements lying in vertex
groups, any reduced sequence given by the theorem \ref{conjthm}
has length at most $k$. It is an easy exercice to verify that one
recovers the original definition of Sela (\cite{acyl}).

\begin{lem}
\label{lgrfini}
 Let $M$ be as above. The JSJ splitting of $\P(M)$ is 4-acylindrical.
\end{lem}

\noindent {\bf Proof.} Let $(\cal{G},X)$ denotes the JSJ splitting
of $\P(M)$. If the splitting is trivial (\textsl{i.e.} $X$ is
reduced to a point), then obviously it is 0-acylindrical and the
conclusion follows, so that we will further suppose that this case
doesn't occur. According to the theorem \ref{conjthm}, there
exists a path in $X$ : $\cal{P}=(v_{\a_0},e_{\b_1},v_{\a_1},\ldots
, e{\b_p},v_{\a_p})$, and a sequence $(c_1,c_2,\ldots , c_p)$ with
$\forall i=1,2,\ldots ,p$, $c_i\in G(e_{\b_i})$, $u\in
G(v_{\a_0})$, $v\in G(v_{\a_p})$, such that $u\sim
c_1^-=\vf_{\b_1}^-(c_1)$ in $G(v_{\a_0})$, $v\sim c_p^+$ in
$G(v_{\a_p})$, and for $i=1,2,\ldots ,p-1$, $c_i^+\sim c_{i+1}^-$
in $G(v_{\a_i})$.  We can also suppose, if for some $i$,
$e_{\b_i}=\bar{e}_{\b_{i+1}}$, that $c_i^+$ and $c_{i+1}^{-}$ are
not conjugate in $G(e_{\b_i})^+=G(e_{\b_{i+1}})^-$ because
otherwise one can  shorten the path and the sequence, while
continuing to verify the above conditions. Hence we can apply
proposition \ref{annuli}, which implies that for $i=1,2,\ldots
p-1$, none of the $M(v_{\a_i})$ is hyperbolic, and hence they are
Seifert fiber spaces. Remark also that  $M(v_{a_i})$ and
$M(v_{a_{i+1}})$ cannot both be twisted $I$-bundles over $\KB_2$,
because otherwise, $M$ would be obtained by gluing two twisted
$I$-bundles over $\KB_2$ along their boundary. Moreover, for
$i=1,2,\ldots , p-2$, $M(v_{\a_i})$ and $M(v_{\a_{i+1}})$ cannot
be both non $I$-twisted bundles over $\KB_2$ : otherwise,
$c_i^+\sim c_{i+1}^-$ in $G(v_{\a_i})$ and $c_{i+1}^+\sim
c_{i+2}^-$ in $G(v_{\a_{i+1}})$ ; hence according to proposition
\ref{annuli}, $c_{i+1}^-$ and $c_{i+1}^+$ must lie respectively in
the fibers of $G(v_{\a_i})$ and $G(v_{\a_{i+1}})$, but this fact
contradicts the lemma \ref{fibers}. Hence, for $i=1,2,\ldots p-1$,
the successive pieces $M(v_{\a_i})$ must be alternatively
$I$-twisted bundles over $\KB_2$, and Seifert pieces which are not
$I$-twisted bundles over $\KB_2$. In fact, since twisted
$I$-bundles have only one boundary component, if for some
$i=1,2,\ldots p-1$, $M(v_{\a_i})$ is a twisted $I$-bundle over
$\KB_2$, then necessarily $e_{\b_i}=\bar{e}_{\b_{i+1}}$, and
$v_{\a_{i-1}}=v_{\a_{i+1}}$ ; hence the pieces
$M(v_{\a_1}),M(v_{\a_2}),\ldots , M(v_{\a_{p-1}})$ are
alternatively a same Seifert piece which is not a twisted
$I$-bundle over $\KB_2$, and possibly several twisted $I$-bundles
over $\KB_2$ each of them being glued to this last non twisted
$I$-bundle piece.

Now suppose $p\geq 5$ ; then without loss of generality, one can
suppose that $v_{\a_1}=v_{\a_3}$,  $M(v_{\a_1})=M(v_{\a_3})$ is a
Seifert piece which is not a twisted $I$-bundle over $\KB_2$,
while $M(v_{\a_2})$ is a twisted $I$-bundle over $\KB_2$. Now,
$c_1^+\sim c_2^-$ in $G(v_{\a_1})$, $c_2^+\sim c_3^-$ in
$G(v_{\a_2})$, and $c_3^+\sim c_4^-$ in $G(v_{\a_1})$. Consider
the canonical presentation $<a,b|aba^{-1}=b^{-1}>$ of
$G(v_{\a_2})$, then according to lemma \ref{annuli}, for some
integers $n,m$, $c_2^+=a^nb^m$ and $c_3^-=a^nb^{-m}$ with $m\not=
0$ (otherwise they would be conjugate in $G(e_{\b_2})^+=<a^2,b>$).
But this last lemma implies also that
$c_2^-=\vf_{\b_2}^{-1}(c_2^+)$ and $c_3^+=\vf_{\b_3}(c_3^-)$ both
lie in the fiber of $G(v_{\a_1})$. This is possible only if
$\vf_{\b_2}$ sends the fiber to the subgroup $<b>$ of
$G(v_{\a_2})$. But this would contradict the lemma \ref{fibers},
since $<b>$ is the fiber of $G(v_{a_2})$ for one of its two
Seifert fibrations. Hence $p\leq 4$ which concludes the
proof.\hfill\bs\\

\subsection{Processing cyclically reduced forms}

Suppose $M$ is a Haken closed manifold whose group $\P(M)$ is
given by its canonical presentation. We have described in \S 3.1,
how one can, given a word $\w$ on the canonical generators of
$\P(M)$, $\w\not= 1$, find a cyclically reduced form whose label is a conjugate of $\w$. This process is
constructive. In order to make use of theorem \ref{conjthm} in a
constructive way, we need to have an algorithmic process to
transform an arbitrary form into a cyclically reduced one. This is
 the aim of this section.

 We claim that in order to perform such a process, it suffices to have a solution to
the generalized
 word problem of $T$ in $\P(N)$,
  for any piece $N$ and
 any boundary subgroup $T$ of $\P(N)$, that is, an algorithm which decides for any $u\in
 \P(N)$ given as a word on the generators of $\P(N)$, whether $u\in T$ or not.
 For, suppose $(\cal{C},\mu)$ is a form with label $\w$. Any
 time $\cal{C}$ contains a sub-path of the form $(v',e,v,\bar{e},v')$,
 check with a solution to the generalized word problem of $G(e)^+$ in $\P(M(v))$ if the
 label $u$ of the vertex $v$
  is an element of $G(e)^+$. Then in this case, replace in
 $\cal{C}$ the subpath $(v',e,v,\bar{e},v')$ with label $(u_1,u,u_2)$ with $(v')$ and label
 $u_1.\vf_e^{-1}(u).u_2$ which is an element of $\P(M(v'))$ (if
 $\cal{C}$ has length 2, replace $(v',e,v,\bar{e})$ with $(v')$ and label $u_1.\vf_e^{-1}(u)$).
 One obtains a shorter cyclic form, with label a word equal in
 $\P(M)$ to $\w$. Perform this process as long as it is possible with $(C,\mu)$ and all of its
 cyclic conjugates.
 Since the length strictly decreases it will necessarily  stop.
  The cyclic form obtained is cyclically reduced. Its
 label is a conjugate in $\P(M)$ of the label $\w$ of $(\cal{C},\mu)$.

A solution to the generalized word problem can be easily found
using the solution to the word problem in the $\P$ of the pieces,
as well as the algebraic properties of boundary subgroups, seen in
the last section.

\begin{prop}
\label{gwp}
 Suppose $N$ is a piece obtained in the non trivial
decomposition of a Haken closed manifold. Suppose $T$ is any
boundary subgroup of $\P(N)$. Then one can effectively decide for
any $u\in \P(N)$, if $u\in T$ or not ; in other words the
generalized word problem of $T$ in $\P(M)$ is solvable.
\end{prop}

\noindent {\bf Proof.} Obviously, we can suppose that $N$ is not
homeomorphic to $S^1\times S^1\times I$. Suppose $N$ is
hyperbolic. Then according to proposition \ref{annuli}, for any
non trivial element $t$ of $T$, the centralizer of $t$ in $\P(N)$
is precisely $T$. Suppose $N$ is a Seifert fibered space, and is
not the twisted $I$-bundle over $\KB_2$. Then, with proposition
\ref{annuli}, for any element $t\in T$ which does not lie in the
fiber $<h>$, the centralizer of $t$ in $\P(N)$ is precisely $T$.
In each  case, to decide if an element $u\in \P(N)$ lies in $T$,
it suffices to apply the solution to the word problem in $\P(M)$
(cf. \cite{wald}) to decide whether  $ut=tu$ or not for such a
$t\in T$.

In the case of the twisted $I$-bundle over $\KB_2$, an element in
$\P(N)=<a,b|aba^{-1}=b^{-1}>$ lies in $\P(\partial{N})$ exactly
when it is written as a word with an even number of occurrence of
the generator $a$ or its inverse, which can be easily
checked.\hfill\bs

As explained above, one immediately obtains the corollary :

\begin{cor}
\label{ccproc}
 Let $M$ be a Haken closed manifold. If $\P(M)$ is given by its canonical presentation,
then one can, given a word $\w$ on the generators, algorithmically
find a cyclically reduced form for $\w$.
\end{cor}

\subsection{The core of the algorithm.} We can now give the main
algorithm to solve the conjugacy problem in the group $\P(M)$ of a
Haken closed manifold $M$. This algorithm uses solutions to the
conjugacy, the boundary parallelism, and the 2-cosets problems in
the groups of each piece, as well as a solution to the word
problem in $\P(M)$. Hence the conjugacy problem
in $\P(M)$ reduces to the conjugacy, boundary parallelism, and
2-cosets problems in the groups of the pieces, which will be
solved latter. We need first to establish two lemmas, to define
correctly the boundary parallelism and 2-cosets problems, and also
to deal with them. The following two lemmas are essential, and are
direct consequences of the algebraic property established in
proposition \ref{annuli}.

In the following  $N$ stands for a piece in the JSJ decomposition
of a Haken closed manifold $M$.

\begin{lem}
\label{ctw}
 For any boundary subgroup $T$ of $\P(N)$, and for
any element $\w\in \P(N)$, define the subset of $T$, $C_{T}(\w)=\{
c\in T \,|\, \w\sim c\ \;\text{in $\P(N)$}\}$.
\begin{itemize}
\item{} If $N$ is hyperbolic, then $C_T(\w)$ is either empty or a singleton.
\item{} If $N$ admits a Seifert fibration then $C_T(\w)$ has
cardinality at most 2.
\end{itemize}
\end{lem}

\noindent {\bf Proof.} Suppose there exists two distinct elements
$c_1$ and $c_2$ in $C_T(\w)$. Then $c_1$ and $c_2$ are conjugate
in $\P(N)$. According to the proposition \ref{annuli}, this cannot
happen if $N$ is hyperbolic, which proves the first assertion. So
$N$ must admit a Seifert fibration. Suppose first that $N$ is not
the twisted $I$-bundle over $\KB_2$. Then necessarily (proposition
\ref{annuli}) $c_1= c_2^{\pm 1}$, and hence $C_T(\w)$ is of
cardinality at most 2. In the case of the twisted $I$-bundle over
$\KB_2$, $c_1=a^{2n}b^p$ and $c_2=a^{2n}b^{\pm q}$ for some
integers $n,p$, and the same conclusion holds.\hfill \bs

\begin{lem}
\label{ctww}
 For any boundary subgroups $T,T'$ of $\P(N)$
(possibly $T=T'$), and for any elements $\w,\w'\in \P(N)$, define
the subset of $T\times T'$, $C_{T,T'}(\w,\w')=\{ (c,c')\in T\times
T'\, |\, \w=c.\w'.c'\}$. Suppose moreover that in  case $T=T'$,
$\w'$ does not lie in $T$.
\begin{itemize}
\item{} If $N$ is hyperbolic,
$C_{T,T'}(\w,\w')$ is either empty or a singleton.
\item{} If $N$ admits a Seifert fibration, and is not the twisted $I$-bundle over $\KB_2$,
let $<h>$ denote the fibre and $C$ the canonical subgroup of
$\P(N)$. Then either $C_{T,T'}(\w,\w')$ is empty or
$C_{T,T'}(\w,\w')=\{ (ch^n,c'h^{-\e.n})\, |\, n\in \Z\}$ with
$\e=\pm 1$ according to whether $\w'\in C$ or not.
\item{} If $N$ is the twisted $I$-bundle over $\KB_2$, then either
$\w\in T$ and $C_{T,T}(\w,\w')$ is empty, or
$C_{T,T'}(\w,\w')=\{(\w{\w'}^{-1}.a^{2n}b^p,a^{-2n}b^p)\, |\,
n,p\in \Z\}$.
\end{itemize}
\end{lem}

\noindent {\bf Remark.} If $T=T'$ and $\w'\in T$, then either
$\w\in T$ and $C_{T,T}(\w,\w')=\{(t.\w{\w'}^{-1},t^{-1})|t\in T\}$
or $C_{T,T}(\w,\w')$ is empty. This is an obvious consequence of
the
fact that $T$ is abelian.\smallskip\\

\noindent {\bf Proof.} Suppose there exist two distinct elements
$(c_1,c_1')$ and $(c_2,c_2')$ in $C_{T,T'}(\w,\w')$. Then
$c_2^{-1}c_1\in T$ is conjugate in $\P(N)$ to $c_2'{c_1'}^{-1}\in
T'$ by ${\w'}^{-1}$. If $N$ is hyperbolic, according to the
proposition \ref{annuli}, necessarily
$c_2^{-1}c_1=c_2'{c_1'}^{-1}=1$, which contradicts the fact that
$(c_1,c_1')$ and $(c_2,c_2')$ are distinct. Thus in the hyperbolic
case, $C_{T,T'}(\w,\w')$ has cardinality at most 1.

If $N$ admits a Seifert fibration and is not the twisted
$I$-bundle over $\KB_2$, then necessarily (proposition
\ref{annuli}) for some integer $n$, $c_2^{-1}c_1=h^n$ and
$c_2'{c_1'}^{-1}=h^{\e.n}$, with $\e=1$ or $-1$, according to
whether $\w'\in C$ or not. Hence, $c_2=c_1.h^{-n}$ and
$c_2'=h^{\e.n}c_1'=c_1'h^{\e.n}$. Reciprocally, if $(c_1,c_1')$
lie in $C_{T,T'}(\w,\w')$, then so does any couple of the form
$(c_1h^n,c_1'h^{-\e.n})$. Hence, if $C_{T,T'}(\w,\w')$ is
non-empty, it must be of the form $\{ (ch^n,c'h^{-\e.n})\, |\,
n\in \Z\}$ for some $(c,c')\in T\times T'$, with $\e=1$ or $-1$,
according to whether $\w'\in C$ or not.

If $N$ is the twisted $I$-bundle over $\KB_2$, then $\P(N)=<a,b\,
|\, aba^{-1}=b^{-1}>$, and $T=<a^2,b>$ is an abelian (normal)
subgroup of index 2. If $w=c\w'c'$ for some $t,t'\in T$, then if
$\w$ lies in $T$, so does $\w'$. Hence if $\w \in T$, such an
equality cannot occur. Since $\w$ and $\w'$ both lie outside $T$
which has index 2, $\w{\w'}^{-1} \in T$, and thus obviously
$(\w{\w'}^{-1},1)\in C_{T,T}(\w,\w')$.

An element of $\P(N)$ lies in $T$, exactly when it can be written
in the form $a^{2n}b^q$ for some integers $n,q$. Since $\w'\not\in
T$, one easily checks that for all integers $n,p$, the equation
$\w'a^{2n}b^p=a^{2n}b^{-p}\w'$ holds. Proceed as above, and
suppose that $C_{T,T}(\w,\w')$ contains two distinct elements
$(c_1,c_1')$ and $(c_2,c_2')$. Necessarily, $\w$ conjugates
$c_2'{c_1'}^{-1}$ into $c_2^{-1}c_1$, and hence with proposition
\ref{annuli}, $c_2'{c_1'}^{-1}=a^{2n}b^p$ and
$c_2^{-1}c_1=a^{2n}b^{-p}$ for some integers $n,p$. So,
$c_2'=a^{2n}b^pc_1'=c_1'a^{2n}b^p$, and $c_2=c_1a^{-2n}b^p$. Then,
the elements of $C_{T,T'}(\w,\w')$ are all those of the form
$(\w{\w'}^{-1}.a^{2n}b^p,a^{-2n}b^p)$ for some integers $n,p$,
which concludes the proof.\hfill\bs\\

Now we consider in the group $\P(N)$ of any piece $N$, two
decision problems : the boundary parallelism problem, and the
2-coset problem.

\noindent {\bf The boundary parallelism problem.} Let $T$ be a
boundary subgroup of $\P(N)$. Construct an algorithm which for any
$\w\in \P(N)$, determines $C_T(\w)$, i.e. find all elements of $T$
conjugate to $\w$ in
$\P(N)$.\smallskip\\
{\bf The 2-coset problem.} Let $T,T'$ be two boundary subgroups of
$\P(N)$ (possibly $T=T'$). Construct an algorithm which for any
couple of elements $\w,\w'\in \P(N)$, determines
$C_{T,T'}(\w,\w')$, i.e.
finds all the couples $(c,c')\in T\times T'$ such that $\w=c.\w'.c'$ in $\P(N)$.\\

We are now able to show that if one can solve the conjugacy,
boundary parallelism and 2-cosets problems  in the groups of all
the pieces obtained in the JSJ decomposition of $M$, then one can
solve the conjugacy problem in $\P(M)$. The cases of $S^1\times
S^1$-bundles over $S^1$, and of two twisted $I$-bundles over
$\KB_2$ glued along  their boundary, are rather easy to deal with,
and a solution to the conjugacy problem in their respective groups
will be sketched in \S 7.

Recall that we suppose that a canonical presentation of $\P(M)$ is
given. Elements of $\P(M)$ are given as words on the canonical
generators.

\begin{thm}
\label{reduce}
 The conjugacy problem in the group of a Haken closed
manifold $M$ which is neither a $S^1\times S^1$-bundle over $S^1$
nor obtained by gluing two twisted $I$-bundles over $\KB_2$ along
their boundary, reduces to conjugacy problems, boundary
parallelism problems, and 2-cosets problems, in the groups of the
pieces obtained. In other words, if one can solve in each of the
groups of the pieces these three problems, then one can solve the
conjugacy problem in $\P(M)$.
\end{thm}

\noindent {\bf Proof.} We will suppose that each piece admits a
solution to these three last problems, and solve the conjugacy
problem in $\P(M)$. Suppose we are given two words $\w$ and $\w'$
on the canonical generators  and want to decide whether or not
$\w$ and $\w'$ are conjugate in $\P(M)$.

First, we use corollary \ref{ccproc} to find cyclically reduced
forms $(\cal{C},\mu)$ and $(\cal{C}',\mu')$ respectively
associated with $\w$ and $\w'$. Without loss of generality we will
suppose that their labels are precisely $\w$ and $\w'$. According
to theorem \ref{conjthm} we can also suppose that the cyclically
reduced forms obtained both have the same length, because
otherwise $\w$,
$\w'$ are definitely not conjugate in $\P(M)$.\\

Suppose first that $(\cal{C},\mu)$ and $(\cal{C}',\mu')$  both
have length 0. This  happens when the paths $\cal{C}$ and
$\cal{C}'$ of $X$ are reduced to points, say $\cal{C}=(v)$,
$\cal{C}'=(v')$, and thus $\w$ and $\w'$ lie in
 the respective vertex groups $G(v)$ and $G(v')$.

 If $v=v'$, apply the solution to the conjugacy problem in $G(v)$ to
 decide whether or not  $\w$ and $\w'$ are conjugate in $G(v)$.
 In the former case, $\w$ and $\w'$ are conjugate in $\P(M)$, but
 in the latter case one cannot conclude yet, and needs to apply the
 general process as described below.

 For any boundary subgroup $T$ of $G(v)$, use the solution in $G(v)$
 to the
 boundary parallelism problem, to find all elements in $T$,
 conjugate to $\w$. According to the lemma \ref{ctw}, one finds at most two such elements. Apply
 the same process with $\w'$ in $G(v')$. One eventually finds
 $c\in G(e)^-\subset G(v)$ and $c'\in G(e')^-\subset G(v')$, two
 respective conjugates of $\w$ and $\w'$. Then apply the same
 process with $\vf_e(c)\in G(e)^+\subset G(\ex(e))$ and $\vf_{e'}(c')\in G(e')^+\subset G(\ex(e'))$,
 and successively
 with all the boundary
 conjugates obtained, to eventually obtain a labelled path from $v$
 to $v'$ as in theorem \ref{conjthm}
 $(ii)$, in which case $\w$ and $\w'$ are conjugate in $\P(M)$.
 Since at each step one  finds at most two boundary
 conjugates, and since according to the lemma \ref{lgrfini} if such a path exists there exists one
 with length at most 4,
 the process must terminate. According to theorem
 \ref{conjthm}, if $\w$ and $\w'$ are not conjugate in some vertex group $G(v)$, and if
 one cannot find such a path, then
 $\w$ and $\w'$ are definitely not conjugate in $\P(M)$.\\

 Suppose now that $(\cal{C},\mu)$ and $(\cal{C}',\mu')$ both have a length greater
 than 0.
 Up to cyclic conjugation of $(\cal{C}',\mu')$ we can suppose that
 $\cal{C}=\cal{C}'=(v_{\s_1},e_{\tau_1},\ldots
,v_{\s_n},e_{\tau_n})$, because otherwise, according to the
theorem \ref{conjthm} $(iii)$, $\w$ and $\w'$ are definitely not
conjugate in $\P(M)$.

First suppose that the path $\cal{C}$ passes through a vertex $v$
whose corresponding piece $M(v)$ is hyperbolic. By possibly
considering cyclic conjugates of $\cal{C}$, one can suppose that
this arises for $M(v_{\s_1})$. According to the theorem
\ref{conjthm} $(iii)$, if $\w$ and $\w'$ are conjugate in $\P(M)$,
then necessarily, there exists $c_n^+\in G(e_{\tau_n})^+$ which
conjugates $\w'$ into $\w$, and moreover there exists also
$c_1^-\in G(e_{\tau_i})^-$ such that
$\mu_1=c_n^+.\mu_1'.(c_1^-)^{-1}$ in $G(v_1)=\P(M(v_1))$. Then,
since $\mu_1$ and $\mu_1'$ are given, using the solution to the
2-cosets problem in $\P(M(v_1))$ one finds at most one couple of
solutions $(c_n^+,c_1^-)$ (cf. lemma \ref{ctww})). Once we know
$c_n^+\in G(e_{\tau_n})^+$, we can use a solution to the word
problem in $\P(M)$ (cf. \cite{wald}), to decide whether or not
$\w=c_n^+.\w'.{(c_n^+)}^{-1}$ in $\P(M)$. In the former case,
obviously $\w\sim\w'$, but in the latter case, before concluding
one needs to apply the same process with all possible cyclic
conjugates $(\cal{C}'',\mu'')$ of $(\cal{C}',\mu')$ such that
$\cal{C}''=\cal{C}'$. Since they are of finite number,  according
to the theorem \ref{conjthm}, one can finally decide whether
$\w\sim\w'$ or not.

Now suppose that the path $\cal{C}$ only passes through vertices
whose corresponding pieces are Seifert fibered spaces. Suppose
first that $\cal{C}$ is of length more than 1 and contains a
subpath of length 1 $(v_1,e,v_2)$ where neither $M(v_1)$ nor
$M(v_2)$ is homeomorphic to the twisted $I$-bundle over $\KB_2$.
Up to cyclic conjugations we can suppose that this condition
arises for the initial subpath $(v_{\s_1},t_{\tau_1},v_{\s_2})$ of
$\cal{C}$. According to the theorem \ref{conjthm}, if $\w$ and
$\w'$ are conjugate in $\P(M)$, then there exist $c_n^+\in
G(e_{\tau_n})^+$, $c_1^-\in G(e_{\tau_1})^-$ and $c_2^-\in
G(e_{\tau_2})^-$, such that,
\begin{gather*}
\w=c_n^+.\w'.(c_n^+)^{-1}\qquad \text{in $\P(M)$}\\
\mu_1=c_n^+.\mu_1'.(c_1^-)^{-1}\qquad \text{in
$\P(M(v_{\s_1}))$}\qquad (1)\\
\mu_2=c_1^+.\mu_2'.(c_2^-)^{-1}\qquad \text{in
$\P(M(v_{\s_2}))$}\qquad (2)
\end{gather*}
We consider (1) and (2) as equations with respective unknowns the
couples $(c_n^+,c_1^-)$ and $(c_1^+,c_2^-)$. We note $S_1$ and
$S_2$ the sets of couples of solutions. Those sets are either
empty or infinite (cf. lemma \ref{ctww}). In the former case $\w$
and $\w'$ are not conjugate. In the latter case we note $C_1^-$
and $C_1^+$ the subsets of $G(e_{\tau_1})$ defined as the
respective images of $S_1$ and $S_2$ under the maps
$(\vf_{\tau_1}^-)^{-1}\circ \pi_2$ and $(\vf_{\tau_1}^+)^{-1}\circ
\pi_1$, where $\pi_1,\pi_2$ stand for the canonical first and
second projections, and
$\vf_{\tau_1}^-:G(e_{\tau_1})\longrightarrow G(e_{\tau_1})^-$,
$\vf_{\tau_1}^+:G(e_{\tau_1})\longrightarrow G(e_{\tau_1})^+$ are
the monomorphisms associated to the edge $e_{\tau_1}$. According
to the lemma \ref{ctww}, $C_1^-$ is a 1-dimensional affine subset
of the $\Z$-module $G(e_{\tau_1})\simeq \ZZ$, with slope
$(\vf_{\tau_1}^-)^{-1}(h_1)$ where $h_1$ is the class of a regular
fiber in $M(v_{\s_1})$, and similarly $C_1^+$ is a 1-dimensional
affine subset of $G(e_{\tau_1})$, with slope
$(\vf_{\tau_1}^+)^{-1}(h_2)$ where $h_2$ is the class of a regular
fiber in $M(v_{\s_2})$. The possible
$c_1=(\vf_{\tau_1}^-)^{-1}(c_1^-)= (\vf_{\tau_1}^+)^{-1}(c_1^+)$
must lie in $C_1^-\cap C_1^+$, and hence, are solutions in $\ZZ$
of a system $(S)$ of two affine equations. The key point is that,
according to lemma \ref{fibers}, $C_1^-$ and $C_1^+$ must have
distinct slopes, and so the system $(S)$ admits at most one
solution --that one can easily determine--. This gives at most one
element $c_1^-$, which according to the lemma \ref{ctww}, allows
the determination of at most one potential element $c_n^+\in
G(e_{\tau_n})^+$ which may conjugate $\w'$ in $\w$ in $\P(M)$. Now
using a solution to the word problem in $\P(M)$, we only need to
check if $\w=c_n^+.\w'.(c_n^+)^{-1}$ in $\P(M)$. If this does not
happen, then apply the same process to all the cyclic conjugates
of $(\cal{C}',\mu')$, whose underlying loops are equal to
$\cal{C}$ (they are of finite number). If one doesn't find in such
a way an element $c_n^+\in G(e_{\tau_n})^+$ which conjugates $\w'$
in $\w$, then, according to the theorem \ref{conjthm} $(iii)$,
$\w$ and $\w'$ are not conjugate in $\P(M)$.

Suppose now, that $\cal{C}$ and $\cal{C}'$ have length one. Then
$\cal{C}=\cal{C}'=(v,e)$ and so the edge $e$ both starts and ends
in $v$. Hence the piece $M(v)$ has at least two boundary
components, and then cannot be homeomorphic to the twisted
$I$-bundle over $\KB_2$. Now, according to the theorem
\ref{conjthm}, $\w=\mu_1.t_e$ and $\w'=\mu_1'.t_e$, and
$\w\sim\w'$ if and only if there exists $c^+\in G(e)^+$ such that
\begin{align*}
\mu_1.t_e &=c^+.\mu_1'.t_e.(c^+)^{-1}\\
        &= c^+.\mu_1'.\vf_e^{-1}((c^+)^{-1}).t_e\\
        \Leftrightarrow \quad \mu_1 &=
        c^+.\mu_1'.\vf_e^{-1}((c^+)^{-1})\quad\text{in $G(v)$}
\end{align*}
Use the solution to 2-coset problem in $G(v)$ to find all couples
$(c,c_1)$ with $c\in G(e)^+, c_1\in G(e)^-$, such that
$\mu_1=c.\mu_1'.c_1$. according to the lemma \ref{ctww}, the set
of solutions $S$ is either empty or $S=\{(d.h^n,d_1.h^{-\e.n})|
n\in \Z\}$, where $h$ is the class of a regular fiber, and $\e=\pm
1$ according to whether $\mu_1'$ (and $\mu_1$) is in the canonical
subgroup $C$ of $\P(M(v))$ or not. If $S=\emptyset$, definitely
$\w$ and $\w'$ are not conjugate. Otherwise, $\w\sim\w'$ if and only if there
exists $(c,c_1)\in S$ such that $\vf_e(c_1)=c^{-1}$, if and only if there
exists $n\in \Z$, such that in $G(e)^+$,
\begin{align*}
\vf_e(d_1.h^{-\e.n})&=h^{-n}.d^{-1}\\
\Leftrightarrow \qquad \vf_e(d_1).d&=h^{-n}.\vf_e(h^{\e.n})\\
 \Leftrightarrow \qquad \vf_e(d_1).d&=(h^{-1}.\vf_e(h^{\e}))^{n}
\end{align*}
Now $G(e)^+\simeq \ZZ$, and once a base is given, one can write in
additive notations, $\vf_e(d_1).d=(a,b)$ and
$h^{-1}.\vf_e(h^{\e})=(\a,\b)$ which are given and do not depend
of $n$. The relation becomes $(a,b)=n.(\a,\b)$, and  hence $\w\sim
\w'$ if and only if the two vectors $(a,b)$ and $(\a,\b)$ of $\ZZ$
are collinear, which can be checked easily. Thus one can decide in
this case whether $\w\sim\w'$ or not.

Now, the only remaining case, is when $\cal{C}$ is of length
greater than one, and such that  for any   sub-path of length 1
$(v,e,v')$ of  $\cal{C}$, either $M(v)$ or $M(v')$ is homeomorphic
to the twisted $I$-bundle over $\KB_2$. Remark that since the
twisted $I$-bundle over $\KB_2$ has only one boundary component,
if the vertex $v$ appearing in $\cal{C}$ is such that $M(v)$ is
homeomorphic to this last manifold, then it must necessarily appear
in a sub-path of $\cal{C}$ of the form $(e,v,\bar{e})$. Moreover
$\cal{C}$  cannot contain a sub-path of the form
$(v_0,e,v_1,\bar{e})$ where both $M(v_0)$ and $M(v_1)$ are twisted
$I$-bundles over $\KB_2$, except when  $M$ is obtained by
gluing two twisted $I$-bundles over $\KB_2$ along their boundary,
which has been excluded.

Suppose first that $\cal{C}$ has length 2. Then up to cyclic
conjugations $\cal{C}=\cal{C}'=(v_0,e,v_1,\bar{e})$, where $M(v_1)$
is homeomorphic to the twisted $I$-bundle, and $M(v_0)$ is not.
Then, $\w=\mu_0.t_e.\mu_1.t_e^{-1}$ and
$\w'=\mu_0'.t_e.\mu_1'.t_e^{-1}$, and according to the theorem
\ref{conjthm} $(iii)$, $\w\sim\w'$ if and only if there exists $c_0,c_1\in
G(e)$, such that
\begin{gather*}
\mu_0=c_0^-.\mu_0'.c_1^-\qquad \text{in $G(v_0)$}\\
\mu_1=c_1^+.\mu_1'.(c_0^+)^{-1}\qquad \text{in $G(v_1)$}
\end{gather*}
Now use the solution to the 2-cosets problem in $G(v_0)$ to find
 the subset $S$ of
$G(e)^-\times G(e)^-$, of all possible $(c_0^-,c_1^-)$ satisfying
the first above equation. According to  lemma \ref{ctww} one obtains
$S=\{(d_0.h^n,d_1.h^{-\e.n})|n\in \Z\}$, where $h$ is the class of
a regular fiber of $M(v_0)$ and $\e=\pm 1$ according to whether
$\mu_0'$ lies in the canonical subgroup $C$ or not. Pick the base
$a^2,b$ of $G(e)^+$, and using additive notations in $G(e)^+$,
note $\vf_e(h)=(p,q)$, then  $c_0^+=\vf_e(c_0^-)=(\a,\b)+n.(p,q)$
and $c_1^+=\vf_e(c_1^-)=(\g,\d)-\e.n.(p,q)$, for $n\in\Z$, and for
some elements $(\a,\b)$ and $(\g,\d)$ that one directly finds from
$S$, $\vf_e$, and the base. Now, $\w\sim\w'$ exactly when
$\mu_1=c_1^+.\mu_1'.(c_0^+)^{-1}$ in $G(v_1)$. According to the
lemma \ref{ctww}, this happens exactly when $c_1^+$ is the image
of $c_0^+$ by the transformation of $\ZZ$ obtained by composing
the
 linear map
which sends $a^2\longrightarrow a^2$, and $b\longrightarrow
b^{-1}$ followed by the translation of vector
$\mu_1{\mu_1'}^{-1}=(\l,\theta)$. Hence $\w\sim\w'$ if and only if
there exists an integer solution $n$ of the equation
$n.(p,-q)+\e.n.(p,q)=(\g-\a-\l,\d+\b-\theta)$, where
$\e,p,q,\a,\b,\g,\d,\l,\theta$ are given, which can be easily
checked.

Now suppose $\cal{C}$ has length greater than 2. Then up to cyclic
conjugation it must contain as initial sub-path
$(v_0,e,v_1,e^{-1},v_0,\ldots)$ where $M(v_1)$ is homeomorphic to
the twisted $I$-bundle over $\KB_2$, while $M(v_0)$ is not.
According to the theorem \ref{conjthm}, $\w\sim\w'$ if and only if there
exists $c_n^+$ in  $G(e_{\tau_n})^+$ which conjugates $\w'$ to
$\w$. Moreover, necessarily, there exist elements $c_1,c_2\in
G(e)$ and $c_3\in G(e_{\tau_3})$, such that
\begin{gather*}
\mu_1=c_n^+.\mu_1'.c_1^-\qquad \text{in $G(v_0)$}\\
\mu_2=c_1^+.\mu_2'.c_2^+\qquad \text{in $G(v_1)$}\\
\mu_3=c_2^-.\mu_3'.c_3^-\qquad \text{in $G(v_0)$}
\end{gather*}
with $c_1^-,c_2^-\in G(e)^-$, $c_3^-\in G(e_{\tau_3})^-$, and
$c_1^+,c_2^+\in G(e)^+$. Note $C_1$ the set of possible $c_1^+\in
G(e)^+$ such that $c_1^-$ verifies the first equation, and $C_2$
the set of possible $c_2^+\in G(e)^+$ such that
$c_2^-=\vf_e^{-1}(c_2^+)$ satisfies the last equation, and use the
solutions to the 2-cosets problem, to find them. We will suppose
that they are both non-empty because otherwise $\w$ and $\w'$ are
not conjugate. Note $h$ the class of a regular fiber in
$\P(M(v_0))$. Pick the base $a^2,b$ of $G(e)^+$, and use additive
notations. Then according to  lemma \ref{ctww}, $C_1$ and $C_2$
are 1-dimensional affine subsets, $C_1=(\a,\b)+\Z.(p,q)$ and
$C_2=(\g,\d)+\Z.(p,q)$ of $G(e)^+$, where $(p,q)$ stands for the
natural image of $h$ under $\vf_e$. But now, necessarily, if
$\w\sim \w'$ then $\mu_2=c_1^+.\mu_2'.c_2^+$ in $G(v_1)$.
According to the lemma \ref{ctww}, this happens exactly when
$c_1^+$ is the image of $c_2^+$ by the transformation of $\ZZ$,
composed of  the linear transformation defined by
$a^2\longrightarrow a^{-2}$ and $b\longrightarrow b$ followed by
the translation of vector $\mu_1{\mu_1'}^{-1}=(\l,\theta)$. Hence,
if $c_1^+=(\a,\b)+n.(p,q)$ and $c_2^+=(\g,\d)+m.(p,q)$ this gives
rise to the equation with the unknowns $n,m\in\Z$,
$$n.(p,q)+m.(p,-q)=(\l-\g-\a,\d+\theta-\b)$$
Now, according to the lemma \ref{fibers}, $(p,q)$ cannot be a
regular fiber of $M(v_1)$, and hence (remember the two Seifert
fibrations of the twisted $I$-bundle over $\KB_2$) neither $p$ nor
$q$ is null. Hence this gives rise to a system of two affine
equations, which admits at most one couple of integer solutions
$(n,m)$. Now, once we know $n$, we know $c_1^-$, and consequently
we know $c_n^+$ (according to the lemma \ref{ctww}). To decide if
$\w\sim \w'$, it suffices to check with the solution to the word
problem in $\P(M)$, whether $\w=c_n^+.\w'.(c_n^+)^{-1}$ or not.

Hence, given $\w$ and $\w'$ in $\P(M)$, by applying this process,
 one can
decide whether $\w\sim\w'$ or not. Hence the conjugacy problem in
$\P(M)$ is solvable, which concludes the argument.\hfill\bs\\

The rest of our work will now consist  in finding solutions to the
boundary parallelism, 2-cosets, and conjugacy problems in the
groups of a Seifert fiber space, or a hyperbolic 3-manifold with
finite volume, as well as in solving the conjugacy problem in the
few remaining cases of $S^1\times S^1$-bundles overs $S^1$ or
manifolds obtained by gluing two twisted $I$-bundles over $\KB_2$
(cf. \S 7).

\section{\bf The case of a Seifert fibered space} This section is
devoted to obtaining  the needed algorithms in the group of a
Seifert fiber space.  We focus essentially on the boundary
parallelism and 2-cosets problem : almost all Seifert fiber spaces
have a biautomatic group, and hence a solvable conjugacy problem ;
the only remaining case --the one of manifolds modelled on $NIL$--
can be treated easily, and will only be sketched in \S 5.3.

\subsection{Preliminaries}
Recall that if $M$ is a Seifert fiber space, any regular fiber
generates a cyclic normal subgroup $N$ called the {\it fiber}.
Moreover $N$ is infinite exactly when $\P(M)$ is infinite (cf. \S
4.1).
$$1\longrightarrow N \longrightarrow \P(M) \longrightarrow \P(M)/N
\longrightarrow 1$$
 Note also that the property of having a group
which contains a normal cyclic subgroup characterizes among all irreducible
3-manifolds with infinite $\pi_1$ those admitting a Seifert fibration (known as the "Seifert fiber space
conjecture", this has been recently solved as the result of a
collective work, including Casson, Gabai, Jungreis, Mess and
Tukia).

The quotient group $\P(M)/N$, is one of a well-known class of
groups, called {\it Fuchsian groups} in the terminology of
\cite{js} (be aware that this definition is "larger" than the
usual definition of a Fuchsian group, as a discrete subgroup of
$PSL(2,\R)$).  If $\r :\P(M)\longrightarrow \P(M)/N$ is the
canonical surjection, and if we note $\ul{u}=\r(u)$, then
$\P(M)/N$ admits one of the following presentations, according to
whether the base of $M$ can be oriented or not :
\begin{gather*}
<\ul{a}_{1},\ul{b}_{1},\ldots,\ul{a}_{g},\ul{b}_{g},\ul{c}_{1},
\ldots,\ul{c}_{q},\ul{d}_{1} ,\ldots,\ul{d}_{p}\mid
\ul{c}_{j}^{\alpha _{j}}=1,\ (\prod_{i=1}^{g} \lbrack
\ul{a}_{i},\ul{b}_{i}\rbrack)\ul{c}_{1}\cdots \ul{c}_{q}
\ul{d}_{1}\cdots \ul{d}_{p}=1>\\
<\ul{a}_{1},\ldots
,\ul{a}_{g},\ul{c}_{1},\ldots,\ul{c}_{q},\ul{d}_{1}
,\ldots,\ul{d}_{p}\mid \ul{c}_{j}^{\alpha _{j}}=1,\
\ul{a}_{1}^{2}\ul{a}_{2}^2\cdots \ul{a}_{g}^{2}.\ul{c}_{1}\cdots
\ul{c}_{q}.\ul{d}_{1}\cdots \ul{d}_{p}=1 >
\end{gather*}
Such groups can be seen as the $\P$ of compact Fuchsian
2-complexes (cf. \cite{js}), or in a more modern terminology, as
$\P^{orb}$ of compact 2-orbifolds whose singular sets consist only
of a finite number of cone points. In this terminology, $M$
inherits a structure of $S^1$-bundle over such an orbifold (cf.
\cite{scott}).

When $M$ has non-empty boundary, the quotient group $\P(M)/N$ is
particularly simple. Indeed, the last relation of the above
presentations, can be transformed into a relation of the form
$\ul{d}_i=\w$ (for some $i=1,2,\ldots ,p$), where $\w$ is a word
which does not involve the letter $\ul{d}_i$ or its inverse ; this
allows the use of Tietze transformations, to discard this relation
together with the letter $\ul{d}_i$. Hence, $\P(M)/N$ is the free
product of the cyclic groups generated by the remaining
generators. The element $\ul{d}_i$ is represented by the word
$\w$, which is cyclically reduced and has length greater than 1
(in the sense of the free product decomposition).

\subsection{Solving the boundary parallelism and 2-cosets
problems}

We solve the boundary parallelism and 2-cosets problems. In both
cases the idea is to reduce to the Fuchsian group $\P(M)/N$, which
easily provides solutions.

\begin{prop}
\label{bound_seif}
 The boundary parallelism problem is solvable in
the group of a Seifert fibered space with non-empty boundary.
\end{prop}

\noindent {\bf Proof.}  We construct an algorithm which solves
this problem. Remark first that in the cases of $S^1\times D^2$,
$S^1\times S^1\times I$, and of the twisted $I$-bundle over
$\KB_2$, the solutions are obvious, so that we can exclude these
cases.

Suppose $T$ is a boundary subgroup of $\P(M)$, generated by
$d_1,h$, and that $u\in \P(M)$ is the conjugate of an element of
$T$, say $u\sim d_1^\a h^\b$ for some integers $\a,\b$. Hence,
$\ul{u}\sim \ul{d}_1^\a$ in $\P(M)/N$. Since, $M\not\simeq
S^1\times D^2$, $M$ is $\partial$-irreducible, and thus $\ul{d}_1$
has infinite order.

Since $M$ has non-empty boundary, $\P(M)/N$ is a free product of
cyclic groups. The element $\ul{d}_1$ is either a canonical
generator, or a cyclically reduced word of length greater than 1.
Now, using the conjugacy theorem in a free product (cf.
\cite{mks}), one can easily determine if $\ul{u}\sim \ul{d}_1^\a$
in $\P(M)/N$, for some integer $\a$, and eventually find
$\ul{a}\in \P(M)/N$ which conjugates $\ul{d}_1^\a$ into $\ul{u}$.
If $\ul{u}$ is not conjugate to $\ul{d}_1^\a$ for some $\a$, then
definitely $u$ is not conjugate in $\P(M)$ to an element of $T$.
Else, if
$$\ul{u}=\ul{a}.\ul{d}_1^\a.\ul{a}^{-1} \qquad \text{in $\P(M)/N$}$$
then once $a\in \r^{-1}(\ul{a})$ has been chosen, since
$\ker\,\r=N$,
$$u=a.d_1^\a.a^{-1}.h^\b=a.d_1^\a h^{\e.\b}.a^{-1}\qquad \text{in $\P(M)$}$$
for some $\b\in \Z$ and $\e=\pm 1$ according to whether $a\in C$
or not. Using the word problem solution in $\P(M)$, one can find
$\b\in \Z$, and thus the element $d_1^\a h^\b\in T$ is conjugate
with $u$ in $\P(M)$. According to the lemma \ref{ctw}, if $\a\not=
0$ or if the base of  $M$ is oriented, then it is the unique
element of $T$ conjugate to $u$. Otherwise, $u$ is conjugate only
to the two elements $h^\b$ and $h^{-\b}$ of $T$. \hfill\bs

\begin{prop}
\label{2coset_seif} The 2-cosets problem is solvable in the group
of a Seifert fibered space with non-empty boundary.
\end{prop}

\noindent {\bf Proof.} In the cases of $S^1\times S^1\times I$ and
$S^1\times D^2$ the solutions are obvious, and in the case of the
twisted $I$-bundle over $\KB_2$, the lemma \ref{ctww} implicitly
provides a solution, so that we can exclude these cases.

Let $T$, be a boundary subgroup of $\P(M)$ ; it is a free abelian
group of rank 2 generated by $d_1,h$. Let $u,v\in \P(M)$, we begin
to determine $C_{T,T}(u,v)=\{(t,t')\in T\times T\,|\, u=t.v.t'\}$.
We first use the proposition \ref{gwp} to decide whether $v\in T$
or not.

In the former case, since $T$ is abelian,
$C_{T,T}(u,v)=\{(uv^{-1}.t,t^{-1})\,|\, t\in T\}$, so that we can
now suppose that $v\not\in T$. Suppose that $u=t.v.t'$, where
$t=d_1^\a h^\b$ and $t'=d_1^\g h^\d$, for some integers
$\a,\b,\g,\d$. then,
$$\ul{u}=\ul{d}_1^\a.\ul{v}.\ul{d}_1^\g\qquad \text{in $\P(M)/N$}$$
Since $M$ is $\partial$-irreducible, $\ul{d}_1$ has infinite order
in $\P(M)/N$. Moreover, no power of $\ul{v}$ lies in $<\ul{d}_1>$,
since, indeed $<\ul{d}_1>$ has a trivial root structure in
$\P(M)/N$, and we have supposed that $v\not\in T$. Hence, since
$\P(M)/N$ is a free product of cyclic groups, one can use the
normal form theorem (cf. \cite{mks}), to find, if any, such a
couple $(\a,\g)$. Thus, since $\ker \r=N$,
$$u=d_1^\a.v.d_1^\g h^\theta\qquad \text{in $\P(M)$}$$
for some $\theta\in \Z$, that one can easily find using the
solution to the word problem in $\P(M)$. Hence we have found an
element of $C_{T,T}(u,v)$, and using the lemma \ref{ctww}, we
determine precisely $C_{T,T}(u,v)$.

Now suppose, we want to determine $C_{T,T'}(u,v)$ for some
distinct boundary subgroups $T,T'$. Suppose $T$, $T'$ are
respectively generated by $d_1,h$, and $d_2,h$. The elements
$\ul{d_1}$ and $\ul{d_2}$ have infinite order in $\P(M)/N$, and
then, using the free product structure, we find, if any, a couple
of integers $(\a,\g)$, such that, $\ul{u}=\ul{d}_1^\a
.\ul{v}.\ul{d}_2^\g$. Then, we can apply the same process as
before to find precisely $C_{T,T'}(u,v)$\hfill\bs

\subsection{Solving the conjugacy problem}

In almost all cases, if $M$ is a Seifert fibre space, $\P(M)$ is
biautomatic, and hence admits a solution to conjugacy problem (cf.
\cite{nr1}, \cite{nr2}) ; the remaining cases are those of
(closed) Seifert fibre spaces modelled on $NIL$ geometry, that is
$S^1$-bundles over a flat orbifold with non zero euler number.
Anyway the conjugacy problem in groups of Seifert fiber spaces can
be easily solved by direct methods ; this is neither difficult nor
surprising, and we will only sketch a proof. The inquiring reader
might refer to \cite{thesis} for a detailed solution.

Conjugacy problem in $\P(M)$ easily reduces to conjugacy problem
in $\P(M)/N$ and to the problem consisting in determining
canonical generators of the centralizer of any element of
$\P(M)/N$ (it can be cyclic, $\ZZ$ or the group of the Klein
bottle). For suppose $u,v\in \P(M)$ are given, and that we want to
decide whether $u\sim v$ or not. We use a solution to the
conjugacy problem in $\P(M)/N$ to decide whether
$\ul{u}\sim\ul{v}$. If $\ul{u}$ and $\ul{v}$ are not conjugate in
$\P(M)/N$, then $u,v$ are not conjugate in $\P(M)$, else there
exists $\ul{a}\in \P(M)/N$ such that
$\ul{u}=\ul{a}\,\ul{v}\,\ul{a}^{-1}$, and if we choose $a\in
\r^{-1}(\ul{a})$, then $u=ava^{-1}h^p$ in $\P(M)$, for some $p\in
\Z$ that one can determine using the solution to the word problem.
Of course, if $p=0$, $u,v$ are conjugate in $\P(M)$, but if
$p\not= 0$ one cannot at this point conclude. To do so, one needs
to determine the canonical generators of the centralizer
$Z(\ul{v})$ of $\ul{v}$ in $\P(M)/N$. Suppose $Z(\ul{v})$ has
generators $x,y$ ; then $vxv^{-1}x^{-1}=h^{n_1}$ and
$vxv^{-1}x^{-1}=h^{n_2}$ in $\P(M)$ for integers $n_1,n_2$ that
one can determine. Then one can easily see that $u$ and $v$ will
be conjugate in $\P(M)$ exactly when $p,n_1,n_2$ satisfy
arithmetic relations, which depend only on the isomorphism class
of $Z(\ul{v})$, as well as on the memberships of $v,x,y$ of the
canonical subgroup of $\P(M)$.

The problem of determining the centralizer of an element of
$\P(M)/N$, and the conjugacy problem in $\P(M)/N$ can be easily
solved, because $\P(M)/N$ is either finite ($M\approx S^1\times
S^2$, or $\mathbb{P}^3\#\mathbb{P}^3$ or modelled on
$\mathbb{S}^3$), word-hyperbolic ($\H^2\times \mathbb{E}^1$,
$\wt{SL}(2,\R)$), or a Bieberbach group ($NIL$, $\mathbb{E}^3$).
Thus the conjugacy problem can be solved in $\P(M)$.

\section{\bf The case of a hyperbolic piece}

In this section we give solutions to the needed decision problems
in the group of a hyperbolic piece. A solution  to the conjugacy
problem is already well-known, according to the following result,
which is a direct implication of the theorem 11.4.1 (geometrically
finite implies biautomatic)  of \cite{epstein}.
\begin{thm} The group of a hyperbolic 3-manifold with finite volume is
biautomatic and hence has solvable conjugacy problem.
\end{thm}

The two remaining decision problems, namely the boundary
parallelism problem and the 2-coset problem, will be solved using
different approaches. The solution to the boundary parallelism
problem will  involve on one hand word-hyperbolic group theory
and on the other Thurston's surgery theorem in the spirit of
Z.Sela, (\cite{sela}), while the 2-coset problem will involve
relatively hyperbolic group theory in the sense of B.Farb
(\cite{farb}).

We first make some  reviews (far from complete) on word hyperbolic
 groups, in order to recall elementary concepts and to fix
notations.
\subsection{Reviews on hyperbolic groups.}
\label{hyp-gps} To a group $G$ with a fixed finite generating set
$X$, one associates the {\it Cayley graph} $\G=\G(G,X)$, which is
a locally finite directed labelled graph, by choosing a vertex
$\ol{g}$ for each element $g\in G$, and for all $g\in G$ and $s\in
X\cup X^{-1}$ an edge with label $s$, going from $\ol{g}$ to
$\ol{g.s}$. To make $\G$ a metric space we assign to each edge the
length 1, and we define the distance between two points to be the
length of the shortest path joining them. Together with this
metric, $\G$ becomes a proper geodesic space. Since vertices of
$\G$ are in 1-1 correspondence with elements of $G$, the group $G$
inherits  a metric $d_G$, called {\it the word metric}. For an
element $\w\in G$, we note
$|\,\w\,|=d_G(1,\w)=d_\G(\ol{1},\ol{\w})$, while the length of a
word $\w$ on the canonical generators, will be noted by
$\lgr(\w)$. Remark that the group $G$ acts on the left naturally
by isometries on its Cayley graph $\G$.

A finite path $\g$ in $\G$, comes equipped with a label which is a
word on the alphabet $X\cup X^{-1}$, naturally obtained by
concatenating the labels of its edges. Given a vertex $v_0$ in
$\G$, finite paths of $\G$ starting from $v_0$ are in one to one
correspondence with words on the generators. We will often make no
distinction between a finite path and its label, as well as
between an element of $G$ and a vertex of $\G$.

A geodesic metric space is said to be $\d${\it -hyperbolic} if
there exists $\d\geq 0$ such that for any geodesic triangle
$(xyz)$, each of its geodesics, for example $[x,y]$, stays in a
$\d$-neighborhood of the union of the two others, $[y,z]\cup
[x,z]$. Given a finite generating set $X$ of $G$, the group $G$ is
said to be $\d$-hyperbolic (resp. hyperbolic) if its Cayley graph
$\G(G,X)$ is $\d$-hyperbolic (resp. $\d$-hyperbolic for some
$\d\geq 0$). It turns out that the property of being hyperbolic,
does not depend of the choice of a finite generating set $X$ of
$G$. If $\G(G,X)$ is $\d$-hyperbolic, then $\G(G,Y)$ is
$\d'$-hyperbolic ; moreover, once we know a set of words on $X$
representing the elements of $Y$, we can easily give a bound on
$\d'$, in terms of $\d$ and of the maximal length of these words
(cf. \cite{cdp}).

Hyperbolic groups, introduced by Gromov (\cite{gro}) to generalize
fundamental groups of closed negatively curved riemanian
manifolds, have been since largely studied and implemented. It
turns out, that they admit very nice algebraic properties, as well
as a particular efficiency in algorithmic processes. For example
they have solvable word and conjugacy problems, which can be
solved respectively in linear and sub-quadratic times. For basic
facts about hyperbolic groups, and
usual background, we refer the reader to the reference books \cite{sho}, \cite{gro}, \cite{gdlh}, \cite{cdp}.\\

\subsection{Solution to the boundary parallelism problem}

To give a solution to the boundary parallelism problem, we will
make use of the word-hyperbolic group theory, and of  Thurston's
hyperbolic surgery theorem (\cite{th}, see also \cite{bp} theorem
E.5.1).

\begin{thm}[Thurston's hyperbolic surgery theorem]
\label{ths}
 let $M$ be a hyperbolic finite volume 3-manifold with
non-empty boundary. Then almost all manifolds obtained by Dehn
filling on $M$ are hyperbolic.
\end{thm}

\noindent {\bf Remarks :} -- Be careful with the sense of "almost
all". It means that if all the surgery coefficients are big
enough, then the manifold obtained is hyperbolic. If $M$ has only
one boundary component, then "almost all" means "all but a finite
number".\\
-- Let $M$ be a hyperbolic 3-manifold with non-empty boundary, and
let $N$ be a closed hyperbolic 3-manifold obtained by Dehn filling
on $M$. Its fundamental group $\P(N)$ is a cocompact (torsion
free) discrete subgroup of $PSL(2,\C)$, and thus is hyperbolic in
the sense of Gromov. Note $\r :\P(M)\twoheadrightarrow \P(N)$, the
canonical epimorphism. Suppose $T$ is a boundary subgroup of
$\P(M)$, that is $T$ is any maximal parabolic subgroup. Then,
since the hyperbolic structure on $N$ extends hyperbolic
structures on the solid tori used in the surgery, the cores of
surgery are geodesics of $N$ and necessarily $\r(T)$ must be
cyclic infinite.\\

The underlying idea (belonging to Z.Sela) is to get two closed
hyperbolic 3-manifolds $N_1$ and $N_2$ by Dehn filling on $M$, and
to use algorithms in $\P(N_1)$ and $\P(N_2)$  to provide a
solution to the boundary parallelism problem in $\P(M)$. Suppose
$\r_1:\P(M)\twoheadrightarrow \P(N_1)$ and
$\r_2:\P(M)\twoheadrightarrow \P(N_2)$ are the canonical
epimorphisms, and suppose one wants to decide for some $\w\in
\P(M)$ and some boundary subgroup $T\subset \P(M)$ , whether $\w$
is conjugate to an element of $T$ or not. Deciding if $\r_i(\w)$
is conjugate in $\P(N_i)$ to an element of $\r_i(T)$ (for
$i=1,2$), will provide a solution in $\P(M)$. We first need to
establish the following lemma.

\begin{prop}
\label{bpgro}
 Let $G$ be a torsion-free $\d$-hyperbolic group, and $H$ a cyclic
 subgroup of $G$. Then, an arbitrary element of $G$ can be
conjugate to at most one element of $H$. Moreover, there exists an
algorithm, which decides for any element of infinite order $\w\in
G$, if $\w$ is conjugate in $G$ to an element of $H$, and finds an
eventual conjugate of $\w$ in $H$.
\end{prop}

\noindent{\bf Proof.} If $H=\{1\}$ the conclusion comes obviously
with a solution to the word problem, so that we will further
suppose that $H$ is non trivial.

 We first
prove the former part of the assumption. Since $h$ has infinite
order, it fixes two distinct points $h^-$ and $h^+$ in the
boundary $\partial \G$ of the cayley  graph. Suppose that $\w$ is
conjugate with two distinct elements of $H=<h>$, say $h^p$ and
$h^q$, then there exists $\a\in G$ such that $\a
.h^p.\a^{-1}=h^q$. Necessarily the action of $\a$ on $\partial\G$
must preserve ${h^-,h^+}$. Hence (\cite{cdp}, prop. 7.1) $\a$ lies
in a finite extension of $H$. In particular, $\exists\,r>0,s>0$,
such that $\a^r=h^s$. Then in one hand
$\a.h^{p.s}.\a^{-1}=(\a.h^p.\a^{-1})^s=h^{q.s}$, and in the other
$\a.h^{p.s}.\a^{-1}=\a.\a^{p.r}.\a^{-1}=\a^{p.r}=h^{p.s}$, which
implies $p=q$.\smallskip\\
\indent In order to prove the latter part of the assumption, we
make use of the stable norm $\|g\|$ of an element $g\in G$ (cf.
\cite{cdp} \S10.6, \cite{gro}), defined as :
$$\|g\|=\lim_{n\rightarrow\infty} \frac{|\,g^n|}{n}$$
the limit exists since $0\leq |\,g^{n+p}|\leq|\,g^n|+|\,g^p|$, and
$\|g\|$ is indeed the infimum of $\{{|\,g^n|}/{n}\,;\,n>0\}$.
It can be seen easily that the stable norm is invariant by
conjugation, that is if $u=a.v.a^{-1}$, $\|u\|=\|v\|$ (remark that
$|\,|\,u^n|-|\,v^n|\,|\leq 2.|\,a\,|$, divide by $n$, and make $n$
go to infinity).

Now suppose that $\w\sim h^t$ for some $t>0$. Considering a
subsequence with indices $t.n$, one has :
$$\|h\|=\lim_{n\rightarrow\infty}\frac{|\,h^{tn}|}{tn}=\frac{1}{t}\,\lim_{n\rightarrow\infty}\frac{|\,h^{tn}|}{n}$$
But since $\w\sim h^t$,
$$\lim_{n\rightarrow\infty}\frac{|\,h^{tn}|}{n}=\|h^t\|=\|\w\|$$
 and hence :
$$\|h\|=\frac{\|\w\|}{t}$$
 Now the key point is that there exists a
computable constant $K>0$, which only depends on $\d$, such that
any element $g$ of infinite order satisfies $\|g\|\geq K$ (cf.
\cite{gro}, remark p.254, and \cite{delz}, prop. 3.1 for a sketch
of a proof). So we finally get :
$$\frac{|\,\w\,|}{t}\geq\frac{\|\w\|}{t}\geq\|h\|\geq K>0$$
which shows that :
$$t\leq\frac{|\,\w\,|}{K}$$

It is now sufficient to use a solution to the word problem to
compute $|\,\w\,|$, and to decide with a solution to the conjugacy
problem if $\w$ is conjugate with $h^t$ for some $t\in\Z$, with modulus
$|\,t\,|<{|\,\w\,|}/{K}$.\hfill\bs\\

We can now give the solution to the boundary parallelism problem
in the group of a finite volume hyperbolic 3-manifold with
non-empty boundary.

\begin{thm}
The boundary conjugacy problem is solvable in the group of a
finite volume hyperbolic 3-manifold with non-empty boundary.
\end{thm}

\noindent {\bf Proof.} Let $M$ be a finite volume hyperbolic
3-manifold with non-empty boundary, and $\cal{T}\subset \partial
M$ a (toro\" \i dal) boundary component. Enumerate all closed
3-manifolds obtained by Dehn filling on $M$ : each one corresponds
(once bases are given) to a couple of coprime integers for each of
the components of $\partial M$. While continuing the enumeration,
process in parallel for each closed manifold $N$ obtained to the
computation of a finite presentation of $\pi_1(N)$ and then apply
the pseudo-algorithm appearing in \cite{papa}. It checks the
hyperbolicity of $\P(N)$, and if so stops, yielding a constant
$\delta$ such that $\P(N)$ is $\delta$-hyperbolic. Pursue these
parallel process until  you have found two groups $\P(N_1)$,
$\P(N_2)$ --obtained by distinct surgery slopes on $\cal{T}$--
which are hyperbolic, and if so, stop all process. According to
Thurston's hyperbolic surgery theorem, the general process will
terminate. Note $T$ the boundary subgroup of $\P(M)$ associated
with $\cal{T}$ and $g_1,g_2\in T$ the respective elements
associated to surgery slopes (up to inverses) on $\cal{T}$ of
$\P(N_1)$, $\P(N_2)$.

Once $\P(N_1)$, $\P(N_2)$ and constants of hyperbolicity are
given, we can apply a process which allows to decide for any
arbitrary element $\w$, if $\w$ is conjugate to an element of $T$,
and find all such conjugate elements. This process is described
below.

The element $g_1\in \P(M)$ is the class of a simple closed curve
on $\cal{T}$, and hence can be completed to form a base $g_1,h_1$
of $T=\ZZ$. Now consider the canonical epimorphism
$\r_1:\P(M)\longrightarrow \P(N_1)$, $\r_1(T)$ is a cyclic
infinite subgroup of $\P(N_1)$ generated by $\r_1(h_1)=h$. Since
$N_1$ is hyperbolic, $\P(N_1)$  has no torsion. Let $\w$ be an
arbitrary element of $\P(M)$. We want to decide if $\w$ is
conjugate to an element of $T$. If $\r_1(\w)$ is non trivial, then
we can use proposition \ref{bpgro} to find at most one element
$h^p$ conjugate to $\r_1(\w)$ in $\P(N_1)$. If $\r_1(\w)=1$, it is
obviously conjugate to an element of $\r_1(T)$. Hence possible
conjugates of $\w$ in $T$ must be of the form $h_1^{p}.g_1^n$, for
some $n\in \Z$, where $p$ is given. Look at the Cayley graph of
$T$ as naturally embedded in the universal cover $\R^2$ of the
torus $\cal{T}$. Then eventual conjugates in $T$ of $\w$ must lie
on the line, with slope $g_1$ crossing $\ol{h_1^p}$. But applying
the same process in $\P(N_2)$, conjugates of $\w$ in $T$ must lie
in the same time, on the line with slope $g_2$ crossing some given
point. Hence, since the two slopes have been chosen distinct they
must be non-collinear, and one easily finds (by resolving a system
of two linear equations) at most one element of $T$ which can be
conjugate with $\w$ in $\P(M)$. Applying the solution to the
conjugacy problem in $\P(M)$, one determines which element of $T$,
if any, is conjugate with
$\w$.\hfill\bs\\

This method can also be applied to solve the 2-coset problem, but
to do so one needs a refinement of Thurston's surgery theorem,
which asserts that there exists a sequence of closed hyperbolic
manifolds converging to $M$ for the geometric topology. But this
case provides  a difficulty, and gives a solution much less
satisfactory, since the surged closed hyperbolic manifolds do
depend of the element $\w\in G$ that one considers (cf.
\cite{thesis}, \S 4.3). A better approach uses relatively
hyperbolic groups theory in the sense of Farb, as seen in the two
following sections. We first recall in the next section elementary
facts upon relatively hyperbolic groups.

\subsection{Reviews on relatively hyperbolic groups.}

This section comes after \S \ref{hyp-gps}. Consider a finitely
generated group $G$, and finitely generated subgroups
$H_1,H_2,\ldots ,H_n$ of $G$. Start from the Cayley graph $\G$ of
$G$, and for each left coset $g.H_i$ add a new vertex $v(g.H_i)$,
as well as  an edge $e(g.h)$ with length $1/2$ from each vertex
$\ol{g.h}$ such that $h\in H_i$,
 to $v(g.H_i)$. Those new vertices and edges will be called {\it special vertices}
 and  {\it special edges}.
  This gives
rise to a new graph $\hat{\G}$, called the {\it coned-off Cayley
graph}, (which does not have to be locally finite), together with
a natural metric which makes $\hat{\G}$ a (non necessarily proper)
geodesic metric space. Note that $\G$ naturally embeds in
$\hat{\G}$, but that this embedding does not (excepted in the
trivial case involving trivial subgroups) preserve lengths.

The group $G$ is said to be $\d${\it-hyperbolic relatively to}
$H_1,H_2,\ldots ,H_n$, if its coned-off Cayley graph $\hat{\G}$ is
a $\d$-hyperbolic geodesic space. It turns out that this
definition does not depend on the choice of a finite generating
set of $G$.

Suppose $X$ is a finite generating set of $G$, and that one knows
for each $H_i$ a finite set of words $S_i=\{y_{i,j}\,|\,j\}$ on
$X$ generating $H_i$. Given a path $w$ in $\G$, there is a usual
way of finding a corresponding path $\hat{w}$ in $\hat{\G}$.
Processing from left to right, one searches in $w$ a maximal
sub-word on the family $S_i$. For each maximal sub-word say $z_i$
on $S_i$, $z_i$ goes from the vertex $\ol{g}$ to $\ol{g.z_i}$,
replace this path with one edge from $\ol{g}$ to the special
vertex $v(gH_i)$, followed by an edge from $v(gH_i)$ to
$\ol{g.z_i}$ (we  make no distinction between a path and its label
in $\G$). Proceed like this, until it is impossible ; obviously the
process will halt. This replacement gives a surjective map
$\G\twoheadrightarrow \hat{\G}$ which from a path $\w$ in $\G$
gives a path that we shall note $\hat{\w}$ in $\hat{\G}$. If
$\hat{\w}$ passes through some special vertex $v(gH_i)$, we say
that $\w$ (or $\hat{\w}$) {\it penetrates} the coset $gH_i$, or
equivalently that $\w$ (or $\hat{\w}$) penetrates the special
vertex $v(gH_i)$.

The path $w$ of $\G$ is said to be a {\it relative geodesic}, if
$\hat{w}$ is a geodesic of $\hat{\G}$. The path $w$ is said to be
a {\it relative quasi-geodesic}, if $\hat{w}$ is a quasi-geodesic.
A path $w$ in $\G$ (or $\hat{w}$ in $\hat{\G}$) is said to be {\it
without  backtracking}, if for every coset $g.H_i$ that $w$
penetrates, $\hat{w}$ does not return to $g.H_i$ after leaving
$g.H_i$. Obviously a relative geodesic is without backtracking.

To proceed efficiently with relative hyperbolic groups, one needs
a more restricted property, the bounded coset penetration property
:\smallskip\\
{\bf Bounded coset penetration property} (or BCP property for
short) : Let $G$ be a group hyperbolic relatively to
$H_1,H_2,\ldots H_n$. Given finite generating sets for
$G,H_1,\ldots ,H_n$, $G$ is said to satisfy the bounded coset
penetration property, if for every $P\geq 1$, there is a constant
$c=c(P)>0$, so that if $u$ and $v$ are relative $P$-quasigeodesics
in $\G$ without backtracking, and with $d_\G(\ol{u},\ol{v})\leq
1$, then
the following conditions hold : \smallskip\\
-- if $u$ penetrates a coset $gH_i$ but $v$ does not penetrate
$gH_i$, then $u$ travels a $\G$-distance of at most $c$ in
$gH_i$.\smallskip\\
-- If both $u$ and $v$ both penetrate a coset $gH_i$, then the
vertices at which $u$ and $v$ first enter $gH_i$ lie a
$\G$-distance of at most $c$ from each other ; and similarly for
the vertices at which
$u$ and $v$ last exit $gH_i$.\\

It turns out that verifying the BCP property does not depend of a
choice of finite generating sets of $G,H_1,H_2,\ldots ,H_n$ (cf
\cite{farb}).

Our motivation for introducing these notions comes from the
following result (theorem 5.1, \cite{farb}).

\begin{thm}
\label{rhg}
 The fundamental group of a complete finite volume
negatively curved riemanian manifold is hyperbolic relatively to
the set of its cusp-subgroups and satisfies the BCP property. In
particular, the same conclusions hold for fundamental groups of
finite volume hyperbolic 3-manifolds relatively to their boundary
subgroups.
\end{thm}

\subsection{Solution to the 2-cosets problem}

We now give a solution to the 2-coset problem in the group of an
 hyperbolic 3-manifold with finite volume. We make use of the fact
 that $\P(M)$ is hyperbolic relatively to its boundary subgroups,
 and satisfies the BCP property (in fact only the last property is
 necessary). The key point is the following result :

 \begin{lem}
Let $G$ be  a hyperbolic group relatively to its subgroups
$H_1,H_2,\ldots ,H_n$, which satisfies the BCP property. Let
$u,v\in G$, $i,j\in\{ 1,2,\ldots ,n\}$ such that if $i=j$, then
$u$ or $v$ does not lie in $H_i$, and suppose that there exist
$c_1\in H_i$, $c_2\in H_j$ such that $u=c_1.v.c_2$ in $G$.

Then, there exists a constant $K$ which only depends on
$\lgr(u),\lgr(v)$ and on constants related to the relatively
hyperbolic structure, such that $c_1$ and $c_2$ have length at
most $K$ for the word metric $d_G$.
 \end{lem}

 \noindent{\bf Proof.} We suppose finite generating sets are given,
and will see the Cayley graph of $\G$ as (non isometrically)
embedded in the coned-off Cayley graph $\hat{\G}$. Consider words
$u$ and $v$ such that $u=c_1.v.c_2^{-1}$ in $G$, for some $c_1\in
H_1$ and $c_2\in H_2$ (eventually $H_1=H_2$). We choose the words
representing $c_1,c_2$ such that they are labels of relative
geodesics ; hence $\hat{c_1}$ is a  path of $\hat{\G}$ of length 1
starting from $\ol{1}$
going through a special edge to the
special vertex $v(H_1)$ and
going back through a special edge to the vertex $\ol{c_1}$,
and similarily
 $\hat{c_2}$ is a  path of length 1
starting from $\ol{u}$ and ending in $\ol{u.c_2}$, which crosses
the special vertex $v(u.H_2)$.  The relation $u=c_1.v.c_2^{-1}$ in
$G$ gives rise to a quadrilateral in $\G$, with vertices
$\ol{1},\ol{c_1},\ol{u},\ol{uc_2}$ and edges labelled with the
words $c_1,v,c_2,u$, such that those with labels $c_1,c_2$ are
relative geodesics.

For
any path $\alpha\subset\hat{\Gamma}$ and any positive integer $t\leq\lg(\alpha)$, we
will note $\alpha(t)$ the vertex of $\alpha$ such that the sub-path of
$\alpha$ from its origin to $\alpha(t)$ has exactly length $t$.
In the following $u$ and $c_1.v$ stand for the paths in $\hat{\G}$ with labels
respectively $u$ and $v$, going respectively from $\ol{1}$ to $\ol{u}$ and from $\ol{c_1}$ to $\ol{c_1v}$.
Note $\w_1$ a relative geodesic path starting from the origin
$u(0)=c_1(0)=\ol{1}$ of $u$ and ending in $c_1.v(1)$, $\w_2$ a relative
geodesic path starting from $u(1)$  and ending in $c_1.v(1)$,
$\w_3$ a relative geodesic path starting from  $u(1)$  and
ending in $c_1.v(2)$, etc... until $u$ or $c_1.v$ does not have any
vertex left, for suppose $p=\lgr(u)\leq \lgr(v)$, then $\w_{2p-1}$
goes from $u(p-1)$ to $c_1.v(p)$, and then we note $\w_{2p}$ the
relative geodesic from $u(p-1)$ to $c_1.v(p+1)$, $\w(2p+1)$ the
relative geodesic going from $u(p-1)$ to $c_1.v(p+2)$, etc...,
until the last vertex $\ol{u.c_2}$ of $c_1.v$. Hence, finally we obtain
$k=\lgr(u)+\lgr(v)+1$ relative geodesics  $c_1$, $\w_1$, $\w_2,\ldots
, \w_{k-2}$ and $c_2$, such that two successive ones have the same
origins and their extremities lying a distance 1, one to the other, or
conversely.

\begin{figure}[ht]
\center{\includegraphics[scale=0.7]{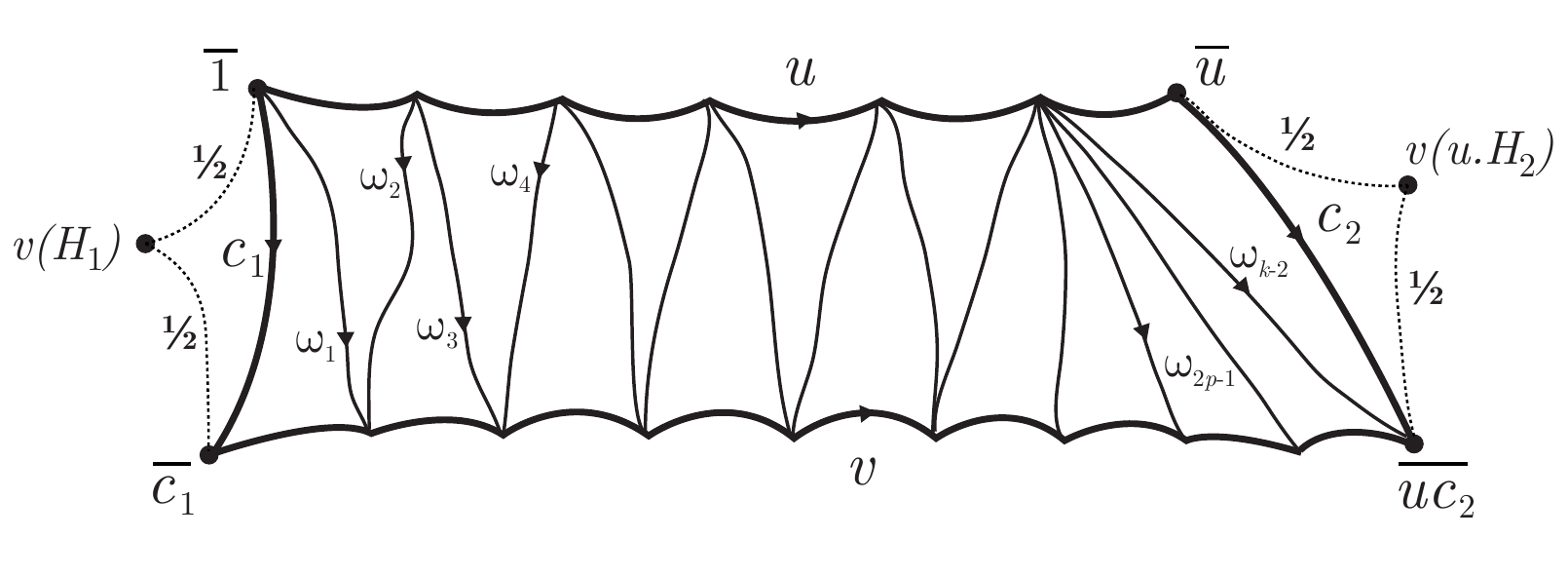}} \caption{}
\end{figure}

Now suppose that $c_1$ has length $L$, then $c_1$ travels a
distance $L$ in the special vertex $v(H_1)$. According to the BCP
property, if $c_1$ has a length greater than $C$, then $\w_1$ must
also cross the special vertex $v(H_1)$, if $c_1$ has length
greater than $3C$, $\w_1$ travels $v(H_1)$ a distance greater than
$C$, and $\w_2$ also crosses $v(H_1)$, etc..., and if $c_1$ has a
length greater than $C.(2(\lgr(u)+\lgr(v))+1)$, $c_2$ must also
cross $v(H_1)$. But since by construction $c_2$ only crosses $v(u.H_2)$ this would
imply $H_1=u.H_2$, and hence $u\in H_1$ and $H_1=H_2$ which
contradicts the hypothesis. Thus we have shown that $|c_1|\leq
C.(2(\lgr(u)+\lgr(v))+1)$ ; the same argument applies to show that
$|c_2|\leq C.(2(\lgr(u)+\lgr(v))+1)$. \hfill\bs

\begin{prop}
The group of a hyperbolic 3-manifold with non-empty boundary and
finite volume has a solvable 2-cosets problem.
\end{prop}

\noindent {\bf Proof.} Using the last lemma, it suffices to
consider all the couples of elements lying in the closed ball of
$G$ with origin $1$ and radius $C.(2(\lgr(u)+\lgr(v))+1)$, and to
use the solution to the word problem to find a possible couple
$(c_1,c_2)$ with $u=c_1.v.c_2$. For such any couple, use the
proposition \ref{gwp} to decide if $c_1\in H_1$ and $c_2\in H_2$,
and conclude.\hfill\bs

\section{\bf Solutions to the conjugacy problem in the remaining
cases}

The remaining cases are :\\
-- A $S^1\times S^1$-bundle over $S^1$. Its group is an HNN
extension of $\ZZ$ with associated isomorphism $\f
:\ZZ\longrightarrow \ZZ$ lying in $SL(2,\Z)$ ; it is indeed the
semi-direct product $\ZZ\rtimes_\f\Z$. Moreover we can suppose
that $\f$ is Anosov (two different irrational eigenvalues), for if
$\f$ is periodic (two complex conjugate eigenvalues, $p$-roots of
the unity) or reducible (one eigenvalue 1 or $-1$), one sees
easily that the
manifold admits a Seifert fibration.\\
-- A manifold obtained by gluing two twisted $I$-bundles over
$\KB_2$ along their (toroidal) boundary. Its group is an
amalgamated product of two copies of $<a,b\,|\,aba^{-1}=b^{-1}>$
along the two copies of the subgroup $<a^2,b>$, with associated
isomorphism lying in $SL(2,\Z)$.

In each case, using the conjugacy theorem in amalgams (cf.
\cite{mks}, \cite{thesis}) or HNN extensions (cf. \cite{ls},
\cite{thesis}), the conjugacy problem reduces to matrix equations
in $SL(2,\Z)$, which can be easily solved.\\

Suppose $M$ is a $S^1\times S^1$-bundle over $S^1$, with an Anosov
associated gluing map. Then $\P(M)=\ZZ\rtimes_\f\Z$ where $\f\in
SL(2,\Z)$ has two distinct irrational eigenvalues. If we note $G$
the first factor of the semi-direct product,  and $t$ a generator
of the second factor, each element of $\P(M)$ can be uniquely
written in the canonical form $u.t^p$, where $u\in G$ and
$p\in\Z$. Now consider two elements $u.t^p$ and $v.t^q$ in
$\P(M)$, and suppose that they are conjugate. Then considering the
homomorphism from $\P(M)$ to $\Z$ which sends $\ZZ$ to $0$ and $t$
to a generator of $\Z$, necessarily $p=q$.

Suppose first that $p=q=0$. Then $u$ and $v$ are conjugate in
$\P(M)$ if and only if there exists $n\in\Z$ such that
$u=\f^n(v)$.  To decide so, consider a base of $G=\ZZ$ constituted
with eigenvectors of $\f$. Then the above equation is equivalent
to the system : $u_1=\l_1^n.v_1$ and $u_2=\l_2^n.v_2$, where
$u=(u_1,u_2)$, $v=(v_1,v_2)$ and $\l_1,\l_2$ are the eigenvalues of
$\f$. The elements $u_1,u_2,v_1,v_2,\l_1,\l_2$ lie in the
extension field $\Q(\sqrt{\Delta})$ where $\Delta$ stands for the
discriminant of the characteristic polynomial of the matrix
associated to $\f$ according to the canonical basis of $\ZZ$. This
system can be easily solved providing a solution in this case.

Suppose now that $p=q$ are distinct from $0$. Using the conjugacy
theorem in an HNN extension, one sees easily that $u.t^p$ and
$v.t^p$ are conjugate if and only if there exists $c\in G$ such
that $u^{-1}v=c.\f^p(c)^{-1}$ up to cyclically conjugating
$v.t^p$. Let us first see how to decide if there exists $c\in G$
such that $u^{-1}v=c.\f^p(c)^{-1}$. Consider the canonical base of
$\ZZ$ ; in this base $u^{-1}v=(n_1,n_2)$, and $\f^p$ has
associated matrix :
$$\mathcal{M}=Mat(\f^p)=\left(\begin{array}{cc}
                        \a & \b \\
                        \g & \d
                        \end{array}
                        \right)$$
with $\a,\b,\g,\d\in \Z$. We look for $c=(x,y)\in\ZZ$ such that
$u^{-1}v=c.\f^p(c)^{-1}$. Then,
$$u^{-1}v=c\f^p(c)^{-1}
        \quad \iff \quad
        \left(
        \begin{array}{c}
        n_1\\
        n_2
        \end{array}
        \right)
                =\left(
                \begin{array}{c}
                x\\
                y
                \end{array}
                \right)
                        -\left(
                        \begin{array}{cc}
                        \a & \b \\
                        \g & \d
                        \end{array}
                        \right)
                                .\left(
                                \begin{array}{c}
                                x\\
                                y
                                \end{array}
                                \right)$$
$$\quad\quad\;\iff \quad\left\lbrace
        \begin{array}{l}
        n_1=(1-\a)x-\b y\\
        n_2= -\g x + (1-\d)y
        \end{array}
        \right.$$
This system has determinant :
$\mathrm{det}(\mathrm{Id}-\mathcal{M})$, which is the value in 1
of the characteristic  polynomial of $\mathcal{M}$. But since $\f$
is Anosov, $\f^p$ is also Anosov and hence does not admit $1$ as
eigenvalue, thus the system admits a unique solution $(x,y)\in
\Q\times \Q$, and if $x$ and $y$ both lie in $\Z$, then $u.t^p$
and $v.t^p$ are conjugate in $\P(M)$. To conclude apply the same
process with all the $p+1$ cyclic conjugates of $v^.t^p$ (they have
the form $\f^q(v).t^p$ for $q=0,1,\ldots ,p$). According to the
conjugacy theorem if one doesn't find in this way that $u.t^p\sim
v.t^p$, then they are definitely not conjugate.
\hfill$\square$\\

Suppose $M$ is obtained by gluing two twisted $I$-bundles over
$S^1$. Note $N_1\approx N_2$ the two $I$-bundles,and $\f$ the
gluing homeomorphism $\f :\partial N_1 \longrightarrow
\partial N_2$. We note
$\P(N_1)=<a_1,b_1|a_1b_1a_1^{-1}=b_1^{-1}>$,
$\P(N_2)=<a_2,b_2|a_2b_2a_2^{-1}=b_2^{-1}>$, for $i=1,2$,
$H_i=<a_i^2,b_i>$, and $\vf : H_1 \longrightarrow H_2$ the
isomorphism induced by $\f$. By fixing respective basis $(a_1)^2,
b_1$ and $(a_2)^2,b_2$ of the free abelian groups of rank two
$H_1$ and $H_2$, $\vf$ can be seen as an element of $SL(2,\Z)$.
With  these bases, note also $\vf_1$ and $\vf_2$ the respective
automorphisms of $H_1$ and $H_2$ associated to the matrix :
$$\left(\begin{array}{cc}
                        1 & 0\\
                        0 & -1
                        \end{array}
                        \right)$$
The group $\P(M)$ is the amalgamated product of $\P(N_1)$ and
$\P(N_2)$ along $\vf$. Each element of $\P(M)$ can be cyclically
reduced in an element which either lies in a factor $\P(N_1)$ or
$\P(N_2)$, or is  of the form : $U=(a_1a_2)^n.u$, with $n\in\Z$
and $u\in H_1=H_2$. We consider $U,V\in\P(M)$, and want to decide
if they are conjugate. According to the conjugacy theorem in
amalgams (cf. \cite{mks}), if $U\sim V$, then  up to cyclic
conjugations, either $U,V$ both lie in a factor, or there exists
$n\in\Z$, and $u,v\in H_1=H_2$, such that $U=(a_1a_2)^nu$ and
$V=(a_1a_2)^nv$.

Suppose first that $U$ and $V$ both lie in a factor, say
$U=a_1^{n_1}b_1^{m_1}$ and $V=a_1^{n_2}b_1^{m_2}$ lie in
$\P(N_1)$, we need to decide if they are conjugate in $\P(N_1)$.
It is an easy exercise to show that $U\sim V$ in $\P(N_1)$ exactly
if either $n_1=n_2$ and $m_1=\pm m_2$, or $n_1=n_2$ is odd and
$m_1=m_2 \mod 2$ ; which can be easily checked. If $U$ and $V$ are
not conjugate in $\P(N_1)$, or if they lie in distinct factors,
then according to the conjugacy theorem, to be conjugate in
$\P(M)$ they must necessarily be conjugate in their respective
factors $\P(N_i)$ to elements of $H_i$ ($i=1$ or 2), and thus,
since $H_i$ is normal in $\P(N_i)$, they must lie in $H_i$. By
eventually composing $U$ or $V$ by $\vf^{-1}$, we will suppose
that both $U$ and $V$ lie in $H_1\subset \P(N_1)$. Applying the
conjugacy theorem in this case one can see easily that $U\sim V$
in $\P(M)$ exactly if there exists an integer $n$, such that (equation $(*)$) :
$(n_2,m_2)=(\vf^{-1}\circ\vf_2\circ\vf\circ\vf_1)^n(n_1,m_1)$.
 Suppose
$$Mat(\vf)=\left(\begin{array}{cc}
                        \a & \b\\
                        \g & \d
                        \end{array}
                        \right)$$
then an easy calculation shows that :
$$\mathcal{M}=Mat(\vf^{-1}\circ\vf_2\circ\vf\circ\vf_1)=\left(\begin{array}{cc}
                        \a\d +\b\g & -2\b\d\\
                        -2\a\g & \a\d +\b\g
                        \end{array}
                        \right)$$
and that the associated endomorphism is either Anosov or reducible
according to whether $\b\g\not=0,-1$ or not. When it is Anosov one
can diagonalise the matrix $\mathcal{M}$, and then easily decide
if a solution $n$ to $(*)$ exists, that is whether $U\sim V$ in
$\P(M)$ or not. When it is reducible, the matrix $\mathcal{M}$ can
only be triangulised, but in this form its diagonal consists only
either of $1$ or of $-1$, and $\mathcal{M}^n$ has a very simple
form which can be easily computed and used to solve $(*)$,
concluding this case.

Suppose now that neither $U$ nor $V$ lie in a factor, and that
they are conjugate in $\P(M)$ ; then for some $p\in\Z$,
$U=(a_1a_2)^p.u$ and $V=(a_1a_2)^p.v$, with $u,v\in H_1$. We note
$\psi= \vf^{-1}\circ\vf_2\circ\vf\circ\vf_1$ ; applying the
conjugacy theorem one obtains that $U\sim V$ in $\P(M)$ if and
only if, up to cyclic conjugation of $V$, there exists $c\in H_1$
such that $vu^{-1}=\psi^p(c).c^{-1}$ (the $p+1$ cyclic conjugates
of $V$ have the form $(a_1a_2)^p\psi^k(v)$ for $k=0,1,\ldots,p$).
This condition is analog to a condition treated above in the case
of a $S^1\times S^1$-bundle over $S^1$, and can be solved in the
same way. There is nevertheless a difference : $\psi^p$ can be
Anosov but also reducible ; as above this last case can be easily
implemented and does not present any difficulty.\hfill$\square$

\vskip 0.4cm

\noindent \scriptsize{{\bf Acknowledgements :} I wish to thank my
advisor Hamish Short for his close attention and his support, as well
as Martin Bridson, Pierre de la Harpe, Gilbert Levitt, Jérome Los, Martin Lustig, for having
accepted to be members of my thesis jury. I want to acknowledge both the referee
of {\it Topology} and Gilbert Levitt, referee of my thesis, for their useful advices.
Finally,
I want to warmly acknowledge Laurence Calon for her help in english and grammar.

During my stay in Geneva in spring-summer 2004, this work has been partially supported by the
'Fond national Suisse de la recherche scientifique'. I wish to thank all the Geneva mathematicians for
their warm welcome.
}

\vskip 0.4cm


\vskip 0.6cm

\noindent Jean-Philippe Pr\'eaux,\smallskip\\
-- Ecole de l'air, CREA, 13661 Salon de Provence air, France.\smallskip\\
-- Laboratoire d'Analyse, Topologie et Probabilités, Centre de math\'ematiques et d'informatique,
39 rue F.Joliot-Curie,
13453 Marseille c\'edex 13,
France.\smallskip\\
\smallskip
E-Mail : \ preaux@cmi.univ-mrs.fr\\
\smallskip
Web : \ \verb/http:///\verb/www.cmi.univ-mrs.fr//\verb$~$\verb/preaux/\\

\end{document}